\DeclareFontFamily{U}{wncy}{}
\DeclareFontShape{U}{wncy}{m}{n}{%
   <5>wncyr5%
   <6>wncyr6%
   <7>wnyr7%
   <8>wncyr8%
   <9>wncyr9%
   <10>wncyr10%
   <11>wncyr10%
   <12>wncyr6%
   <14>wncyr7%
   <17>wncyr8%
   <20>wncyr10%
   <25>wncyr10}{}

\documentclass[a4paper,11pt]{article}

\usepackage{latexsym,graphics,color}
\usepackage{amsmath,amsfonts,amssymb,amsthm}
\newtheorem{thm}{Theorem}[section]
\newtheorem{lem}[thm]{Lemma}

\newtheorem{cor}[thm]{Corollary}
\newtheorem{prop}[thm]{Proposition}

\newtheorem{rem}[thm]{Remark}

\input epsf

\title{On geodesics of phyllotaxis}
\author{Roland Bacher}

\begin{document}
\maketitle
%\par sunflow1.tex dans recherche/sunflowermap 

{\sl Abstract\footnote{
Math. class: 11H31, 52C15 Primary: 92B99}: 
Seeds of sunflowers are often modelled by the map $n\longmapsto 
\varphi_\theta(n)=\sqrt{n}e^{2i\pi n\theta}$ leading to a roughly
uniform repartition with two consecutive seeds separated by 
the divergence angle $2\pi\theta$  for $\theta$ the golden ratio.
We associate to an arbitrary real divergence angle $2\pi \theta$ 
a geodesic path $\gamma_\theta:
\mathbb R_{>0}\longrightarrow \mathrm{PSL}_2(\mathbb Z)\backslash \mathbb H$ 
of the modular curve and use it for local descriptions of the image
$\varphi_\theta(\mathbb N)$ of the phyllotactic map
$\varphi_\theta$.}
\vskip.5cm
Given a real parameter $\theta$, we call the map $\varphi_\theta:\mathbb N
\longrightarrow \mathbb C$ defined by 
$$\varphi_\theta(n)=\sqrt{n}e^{2i\pi\theta n}$$
the \emph{phyllotactic map of divergence angle $2\pi\theta$} 
(measured in radians). The image 
$\varphi_\theta(\mathbb N)$
of a phyllotactic map is the \emph{phyllotactic set} (of parameter $\theta$
or divergence angle $2\pi\theta$).
A phyllotactic set $\varphi_\theta(\mathbb N)$
is uniformly discrete (i.e. two distinct elements of $\varphi_\theta(\mathbb N)$
are at distance at least $\epsilon$ for some strictly positive $\epsilon$)
with uniform density if 
$$\theta=[a_0;a_1,a_2,\dots]=a_0+\frac{1}{a_1+\frac{1}{a_2+\dots}}$$ 
is irrational with bounded coefficients $a_0,a_1,a_2,\dots$
in its continued fraction expansion.

Among all possible parameters, the value given by the golden ratio
$\frac{1+\sqrt{5}}{2}=[1;1,1,1,\dots]$ (or closely related numbers) 
stands out and gives a particularly nice configuration. Figure 1 displays
a few hundred small points of $\varphi_{(1+\sqrt{5})/2}(\mathbb N)$.
%\eject
\begin{figure}[h]\label{figure1}
\epsfysize=8cm
\center{\epsfbox{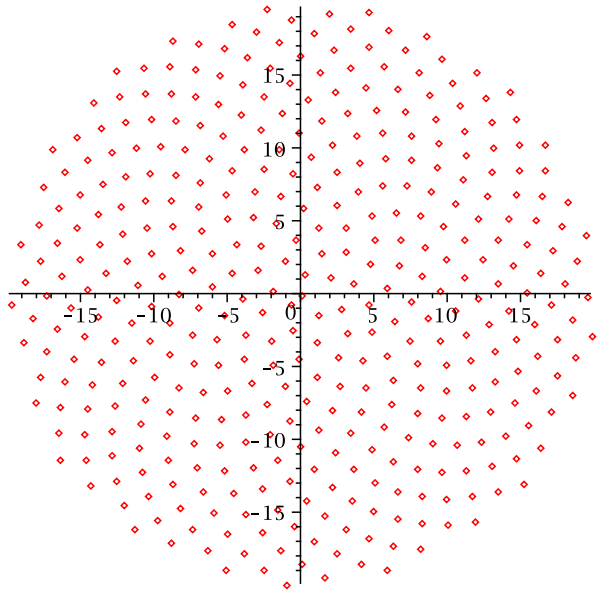}}
\caption{All points of $\varphi_{(1+\sqrt{5})/2}(\mathbb N)$ in the disc 
$\{z\in\mathbb C\mid \vert z\vert\leq 20\}$.}
\end{figure}
%\eject

Finite approximations of $\varphi_{(1+\sqrt{5})/2}(\mathbb N)$
can be observed in capitula (heads) of sunflowers or daisies
(the map $\varphi_{(1+\sqrt{5})/2}$, 
sometimes also called the sun\-flower-map, has 
been proposed in \cite{V} as a model for heads of sunflowers).
Joining close points of $\varphi_{(1+\sqrt{5})/2}(\mathbb N)$
we get \emph{parastichy spirals}
appearing in pairs of crisscrossing families enumerated by two 
consecutive elements of the Fibonacci sequence $1,2,3,5,8,13,21,\dots$. 
Explaining 
the occurence of the golden ratio and of Fibonacci numbers in Botanics
is the goal of Phyllotaxis, see for example Chapter XIV of \cite{DaT}
or \cite{Je} for more recent developments. 
The aim of this paper is to describe
an elegant framework involving hyperbolic geometry. 
The emergence of the golden ratio $\frac{1+\sqrt{5}}{2}$ (or of closely
related numbers) and 
of Fibonacci numbers enumerating families of
parastichy spirals is then a consequence of natural
constraints. 

Reasons for Phyllotaxis should be separated from the mechanisms which
are involved. How Phyllotaxis works is surely best adressed by 
biologists, biochimists or biophysicists.
The reason 
for Phyllotaxis is efficiency of some sort (a precise definition is perhaps
not so easy) which can take several forms.
It is perhaps a physical notion like
energetic efficiency or it involves geometric
quantities like isoperimetry (which leads probably ultimately also
to some kind of energetic efficiency). The link between the two aspects
is natural selection. To say it in a nutshell, ubiquity of phyllotaxis 
involves mathematics: A few
geometric configurations optimize some natural quantities. Thus they
are favoured by living organisms through natural selection. 

Interestingly, the two aspects are spatially
separated: Reasons for Phyllotaxy, 
due to a globally optimized quantity, are of an 
asymptotical nature.
They are best adressed by studying the large part of a plant which
is relatively far from the center consisting of the bud which is 
responsible for the growth-process of a flower.
Asymptotic arguments are thus not a weakness but
are relevant when trying to answer why phyllotaxis occurs (assuming 
Darwin's theory of evolution).

An outline of the paper is as follows:
We associate to a phyllotactic map $\varphi_\theta$ with real parameter 
$\theta\in [0,1)$ the curve
$$(0,+\infty)\ni t\longmapsto \gamma_\theta(t)=
\frac{4i\pi t}{4i\pi \theta t+1}=\frac{16\pi^2\theta t^2+4i\pi (1-\theta)t}
{16\pi^2\theta^2t^2+1}$$
of the Poincar\'e halfplane $\mathbb H=\{z\in\mathbb C\mid \Im(z)>0\}$.
Since $\varphi_\theta=\varphi_{\theta+n}$ for all $n\in\mathbb Z$, we
extend the definition of $\gamma_\theta$ to $\theta\in\mathbb R$ by
setting
$$\gamma_\theta(t)=
\frac{4i\pi t}{4i\pi \lbrace \theta\rbrace t+1}$$
where $\theta=\lfloor \theta\rfloor +
\lbrace \theta\rbrace$ with $\lfloor \theta\rfloor \in\mathbb Z$ and
$\lbrace \theta\rbrace\in[0,1)$
denoting the integral and fractional part of $\theta$. 
It is easy to check that $\gamma_\theta$
defines a geodesic with respect to the usual hyperbolic metric of 
$\mathbb H$. We call $\gamma_\theta$ the 
{\em phyllotactic geodesic} associated to $\varphi_\theta$.
For large $N\in\mathbb N$, the projection of $\gamma_\theta(N)$ 
(often identified with $\gamma_\theta(N)$ in the sequel)
onto the modular curve $\mathrm{PSL}_2(\mathbb Z)\backslash\mathbb H$
classifying complex lattices up to similarity describes, 
up to an affine orientation-preserving 
similarity, the affine lattice obtained by \lq\lq linearizing'' the 
phyllotactic map $\varphi_\theta$ in a neighbourhood of $\varphi_\theta(N)$.
Linearization is an asymptotical construction involving an error of
order $O\left(\frac{1}{\sqrt{N}}\right)$ in a neighbourhood of 
a point $\varphi_\theta(N)$.
In particular, it breaks down for very small values of $N$. 
This failure should have no serious consequences: Indeed, 
the most interesting phyllotactic sets have obviously already good 
packing properties at their center. Moreover, the center 
yields a very small
contribution to interesting quantities like energy, mean isperimetric
values of Voronoi cells etc..

The apparition of the golden ratio $\frac{1+\sqrt{5}}{2}$ 
can now be explained by the fact that the corresponding 
phyllotactic geodesic $t\longmapsto \frac{4i\pi t}{2i\pi(-1+\sqrt{5})t+1}$
avoids the cusp of the modular curve
$\mathrm{PSL}_2(\mathbb Z)\backslash\mathbb H$. This ensures small excentricity
(or, equivalently, 
good isoperimetric constants) for Voronoi cells of the phyllotactic 
set $\varphi_{(1+\sqrt{5})/2}(\mathbb N)$
where the Voronoi cell $V_v$ of a point $v\in S$ with respect to a 
discrete set $S$ in a metric space $E$ is the subset $V_v=\{x\in E\mid
d(x,v)=\min_{y\in S}d(x,y)\}$ of all points closest to $v$. Plants having
seed-areas with good isoperimetric constants should be favoured 
by natural selection since they need less material for constructing 
seed coats around seeds of given volume.
Thus, Diophantine properties of the golden ratio $\tau=\frac{1+\sqrt 5}{2}$
ensure that $\varphi_\tau(\mathbb N)$
(or $\varphi_\theta(\mathbb N)$ for $\theta$ a close relative of 
the golden ratio $\tau$ giving rise e.g. to the sequence
$1,3,4,7,11,18,\dots$ also observed in Phyllotaxis) 
satisfies natural constraints in the class of 
all sets of the form $\varphi_\theta(\mathbb N)$.

We can either consider that the points of the sequence $\varphi_{\theta}(0),
\varphi_\theta(1),\dots$  appear sequentially with constant 
divergence angle $2\pi\theta$ between consecutive points on 
the so-called \emph{ontogenetic spiral}
$t\longmapsto \sqrt{t}e^{2i\pi\theta t}$
or, sticking perhaps closer to biological reality, 
we can consider a sequence of paths $[-n,0]\ni t\longmapsto 
s_{n,\theta}(t)=\sqrt{n+t}e^{2i\pi n\theta }$ describing a situation 
where the $n-$th point has appeared at the origin at time $-n$ in the past.
It has then slowly moved outwards on the halfray $\mathbb R_{\geq 0}e^{2i\pi n\theta}$ until reaching its present location $\varphi_\theta(n)$ at the time $t=0$.
 
The factor $\sqrt{n}$ in the formula for $\varphi_\theta(n)$ 
ensures that there are
$R^2$ points (of roughly equal \lq\lq importance\rq\rq,
measured for example by the area of Voronoi cells) in a disc of large 
radius $R$. Are{\ae} of Voronoi cells defined by 
$\varphi_\theta(\mathbb N)$ are asymptotically equal to 
$\pi$ if $\theta$ is irrational.

Tools and part of the results of this paper can be adapted to a slightly more
general situation given by functions
$t\longmapsto \rho(t)e^{2i\pi t\alpha(t)}$
where $\rho(t)$ is a suitable increasing function and where the local 
divergence angle $2\pi \alpha(t)$ is allowed to vary very slowly.

The sequel of the paper is organized as follows:

Section \ref{sectmain} states the main result.

Section \ref{sectreminders} recalls a few well-known and useful facts
concerning complex lattices and hyperbolic geometry.

Section \ref{sectcontfractions} contains identities 
involving continuous fraction expansions.

Linearizations of phyllotactic sets are
described in Section \ref{sectlinearization}.

We construct the phyllotactic geodesic $\gamma_\theta$ in
Section \ref{sectphyllgeod}. This leads to a proof of Theorem \ref{mainthm}.

Section \ref{sectgeomconstr} describes a construction of a 
slightly different geodesic.

Section \ref{sectmetrprop}
is devoted to metric properties of phyllotactic sets.

Parastichy spirals are defined and studied in Section \ref{sectparast}.

Section \ref{sectcombvor} describes a few combinatorial aspects of 
Voronoi diagrams for phyllotactic sets.

Section \ref{sectchrom} discusses some chromatic aspects
related to local canonical four-colourings of the Voronoi
cells defined by phyllotactic sets.

Section \ref{sectothermod} reviews briefly a few other models
appearing in the literature.

Finally, Section \ref{secttest} discusses a possible experimental 
verification (or refutation) of the existence of a phyllotactic 
geodesic in real sunflower-capitula.

%%%%%%%%%%%%%%%%%%%%%%%%%%%%%%%%%%%%%%%%%%%%%%%%%%%%%%%%%%%%%%%%%%%%%%%%%%

\section{Main result}\label{sectmain}

Let $(E,\hbox{dist})$ be a metric space, $x$ an element of $E$ 
and $\epsilon, R$ two strictly positive real numbers.
Two discrete subsets $A,B$ of $E$ are 
\emph{$\epsilon-$close in the open ball of radius $R$ and center $x$}
if there exists a map $\psi:A'\longrightarrow B'$
which is one-to-one and onto 
between subsets $A'\subset A$ and $B'\subset B$
containing all points in $A$, respectively $B$, at distance at most $R$ 
from $x$ 
and which moves all points of $A'$ by less than $\epsilon$, i.e., 
we have $\hbox{dist}(a,\psi(a))< \epsilon$ for all $a\in A'$.
Intuitively, two discrete sets $A,B$ are $\epsilon-$close
in the open ball of radius $R$ centered at $x$ 
if $A$ and $B$ are \lq\lq equal up to an error
of $\epsilon$'' in (a neighbourhood of) the ball of radius 
$R$ centered at $x$.

\begin{thm}\label{mainthm} 
Given $\epsilon>0$ and $R>0$, there exists an integer $N=N(\epsilon,R)$
such that for every $\theta\in[0,1)$ and for every $n\geq N$, 
the set $\varphi_\theta(\mathbb N)$ is $\epsilon-$close
in the open disc of radius $R$ centered at $\varphi_\theta(n)$ to an affine
lattice in the equivalence class (i.e. up to
orientation-preserving affine similarities) of 
$\mathbb Z+\mathbb Z\frac{4i\pi n}{1+4i\pi\lbrace \theta\rbrace n}$.
\end{thm}

The map
$$(0,\infty)\ni t\longmapsto \gamma_\theta(t)=
\frac{4i\pi t}{1+4i\pi\lbrace \theta\rbrace t}=\frac{16\pi^2\lbrace\theta\rbrace t^2+4i\pi t}{1+(4\pi\lbrace\theta\rbrace t)^2}$$
defines a geodesic of the hyperbolic Poincar\'e halfplane (endowed with the 
hyperbolic metric $\frac{ds}{y}$ at $z=x+iy\in \mathbb H$), 
see for example Lemma \ref{lemgeod}. We call $\gamma_\theta$
the \emph{phyllotactic geodesic} of $\varphi_\theta$.
Notice that the curve $t\longmapsto\gamma_\theta(t)$ has (hyperbolic) speed 
$\frac{1}{\Im(\gamma_\theta(t))}\vert \gamma'_\theta(t)
\vert=\frac{1}{t}$ inversely proportional to $t$.

\begin{rem}\label{remareaVor} (i) Theorem \ref{mainthm} holds for rational $\theta$: 
In this case the phyllotactic geodesic ends up in the cusp of the modular
curve $\mathrm{PSL}_2(\mathbb Z)\backslash\mathbb H$. The associated
affine lattices $\Lambda_{\theta,n}$ degenerate into discrete subgroups of 
rank $1$ in the sense that they intersect a ball of fixed radius $R$ 
centered at an affine lattice point along a translated copy 
of a discrete subgroup having rank $1$.

(ii) Denoting by $V(n)$ the Voronoi cell of $\varphi_\theta(n)\in
\varphi_\theta(\mathbb N)$ we have $\lim_{n\rightarrow\infty}
\mathrm{vol}
V(n)=\pi$ if $\theta$ is irrational.

For rational $\theta=\frac{p}{q}$ with with $q\geq 3$ and 
$p,q$ coprime integers, the union of
Voronoi cells (defined by $\varphi_\theta(\mathbb N)$)
of all points at distance $\leq R$ from the origin
is essentially a regular polygon with $q$ sides and inradius $R$.
This implies $\lim_{n\rightarrow\infty}\mathrm{vol}V(n)=q
\tan\frac{\pi}{q}=\pi+\frac{\pi^3}{3q^2}+O\left(\frac{1}{q^4}\right)$.

For $\theta\in\frac{1}{2}\mathbb Z$, all Voronoi cells 
are unbounded and thus of infinite volume.
\end{rem}

%Since $\mathrm{PSL}_2(\mathbb Z)$ acts isometrically on $\mathbb H$
%(endowed with the hyperbolic metric), 
The existence of phyllotactic geodesics gives a 
measure of \lq\lq similarity'' of phyllotactic sets
in neighbourhoods of $\varphi_\theta(n)$ and $\varphi_{\theta'}(m)$
by considering the hyperbolic 
distance 
$d_{\mathbb H}(\gamma_\theta(n),\mathrm{PSL}_2(\mathbb Z)\gamma_{\theta'}(m))$
between the two orbits $\mathrm{PSL}_2(\mathbb Z)\gamma_{\theta}(n)$ and
$\mathrm{PSL}_2(\mathbb Z)\gamma_{\theta'}(m)$. Statement (ii) of Remark 
\ref{remareaVor} ensures that a small distance 
$d_{\mathbb H}(\gamma_\theta(n),\mathrm{PSL}_2(\mathbb Z)\gamma_{\theta'}(m))$
implies the existence of bijections between 
$\varphi_\theta(\mathbb N)$ and $\varphi_{\theta'}(\mathbb N)$
which are almost isometries in neighbourhoods of  
$\varphi_{\theta}(n)$  and $\varphi_{\theta'}(m)$ if 
$\theta$ and $\theta'$ are irrational and have continued fraction-expansions
with bounded coefficients. More generally, this holds 
if $\gamma_\theta(n)$ (and thus also $\gamma_{\theta'}(m)$)
is far from the cusp of $\mathrm{PSL}_2(\mathbb Z)\backslash \mathbb H$.

%%%%%%%%%%%%%%%%%%%%%%%%%%%%%%%%%%%%%%%%%%%%%%%%%%%%%%%%%%%%%%%%%%%%%%%%%%

\section{Complex lattices and hyperbolic geometry}\label{sectreminders}

For the convenience of the reader, we recall a few elementary and well-known
facts first of the theory of lattices, 
following closely parts of Section 2.2 in Chapter VII of \cite{Serre},
then of hyperbolic geometry, see for example \cite{And}.

\subsection{Lattices of $\mathbb C$}
A \emph{lattice} in $\mathbb C$ is a free additive subgroup generated
by two $\mathbb R-$linearly independent elements $\omega_1,\omega_2$ 
of $\mathbb C$. In the sequel, we consider lattices
only up to orientation-preserving similarities.
Two lattices $\Gamma$ and $\Lambda$ of $\mathbb C$ are thus 
equivalent if $\Lambda=\lambda\Gamma$ for some non-zero constant
$\lambda\in \mathbb C^*$. 
Given a basis $\omega_1,\omega_2$ of a lattice
$\Gamma=\mathbb Z\omega_1+\mathbb Z\omega_2$, 
we consider $z=\frac{\omega_1}{\omega_2}$. Up to replacing,
say, $\omega_1$ by $-\omega_1$, we can suppose that the imaginary 
part $y=\Im(z)$ of $z=x+iy$ is strictly positive.

Thus, a lattice $\mathbb Z\omega_1+\mathbb Z\omega_2=
\left(\mathbb Z+\mathbb Z\frac{\omega_1}{\omega_2}
\right)\omega_2$ is equivalent to the lattice $\Gamma(z)=
\mathbb Z+\mathbb Zz$ generated by $1$ and 
by the element $z=\frac{\omega_1}{\omega_2}$ of the open
upper half-plane $\mathbb H=\{z\in \mathbb C\mid
\Im(z)>0\}$.

Given an unimodular integral matrix $g=\left(\begin{array}{cc}a&b\\c&d
\end{array}\right)\in \mathrm{SL}_2(\mathbb Z)$, the quotient 
$$z'=
\frac{\omega_1'}{\omega_2'}=\frac{a\omega_1+b\omega_2}{c\omega_1+d\omega_2}
=\frac{a\frac{\omega_1}{\omega_2}+b}{c\frac{\omega_1}{\omega_2}+d}$$ 
associated to the basis 
$\omega_1'=a\omega_1+b\omega_2,\ \omega_2'=c\omega_1+d\omega_2$
of a lattice $\mathbb Z\omega_1+\mathbb Z\omega_2$ is obtained from
$z=\frac{\omega_1}{\omega_2}$ by the usual action 
$g.z=\frac{az+b}{cz+d}$ of the modular group $\mathrm{PSL}_2(\mathbb Z)
=\mathrm{SL}_2(\mathbb Z)/\pm \mathrm{Id}$
on $\mathbb H$.

Hence the map $\mathbb Z\omega_1+\mathbb Z\omega_2\longmapsto z=\frac{\omega_1}
{\omega_2}\in \mathbb H$ induces a one-to-one correspondence 
between equivalence classes $\mathbb C^*\Gamma$ of lattices 
and points of the modular curve $\mathrm{PSL}_2(\mathbb Z)\backslash
\mathbb H$,
see Chapter VII, Proposition 2 and Proposition 3 of \cite{Serre}.
A fundamental domain for the action of the modular group 
$\mathrm{PSL}_2(\mathbb Z)$ on $\mathbb H$ is given by the \emph{fundamental
domain}
\begin{eqnarray}\label{deffunddomM}
\mathcal M=\left\lbrace z\in \mathbb C\mid \vert z\vert\geq 1\hbox{ and }
\vert \Re(z)\vert\leq \frac{1}{2}\right\rbrace
\end{eqnarray}
for $\mathrm{PSL}_2(\mathbb Z)\backslash\mathbb H$.
Two elements $z_1,z_2$ of $\mathcal M$ represent the same equivalence-class
of lattices if and only if either $z_2=-\frac{1}{z_1}$ or $z_2=z_1\pm 1$.
The \emph{modular curve} 
$\mathrm{PSL}_2(\mathbb Z)\backslash\mathbb H=\mathcal M/\sim$ 
is a complex orbifold with two 
conical points represented by $\frac{-1+i\sqrt{3}}{2}$ 
(of angle $\frac{2\pi}{3}$
and corresponding to regular hexagonal lattices) and by
$i$ (of angle $\pi$ and corresponding to square lattices)
and with a cusp (corresponding to a neighbourhood of the degenerate case of 
an additive 
subgroup of rank $1$ in $\mathbb C$). 

An  \emph{affine lattice} is a coset $\alpha+\Gamma$ obtained 
by translating a complex lattice $\Gamma\subset\mathbb C$ 
by some vector $\alpha\in \mathbb C$.
We consider affine lattices only up to orientation-preserving affine 
similarities. Equivalence classes of affine lattices 
are also in one-to-one correspondence with elements of the 
modular curve $\mathrm{PSL}_2(\mathbb Z)\backslash\mathbb H$.

%%%%%%%%%%%%%%%%%%%%%%%%%%%%%%%%%%%%%%%%%%%%%%%%%%%%%%%%%%%%%%%%%%%%%%%%%%

\subsection{Hyperbolic geometry on
the Poincar\'e half-plane}

We recall a few facts concerning the hyperbolic Poincar\'e 
half-plane $\mathbb H$, see \cite{And} for an elementary 
introduction to hyperbolic geometry.

The upper half-plane $\mathbb H=\{z\in \mathbb C\ \vert \Im(z)>0\}$
can be turned into a real
hyperbolic simply connected Riemannian manifold of dimension $2$ and of 
constant curvature $-1$ by equipping it with the Riemannian metric 
$(ds)^2=\frac{dx^2+dy^2}{y^2}$ at a point $z=x+iy\in \mathbb H$.
The \emph{Poincar\'e half-plane} is the hyperbolic manifold (still denoted by)
$\mathbb H$ obtained in this way.

The group of all orientation-preserving isometries of the 
Poincar\'e half-plane is given by the set of all M\"obius transformations
$$z\longmapsto\left(\begin{array}{cc}a&b\\c&d\end{array}\right)z=
\frac{az+b}{cz+d}$$
defined by matrices $\left(\begin{array}{cc}a&b\\c&d\end{array}\right)$
in $\mathrm{SL}_2(\mathbb R)$, respresenting elements
in $\mathrm{PSL}_2(\mathbb R)=\mathrm{SL}_2(\mathbb R)/\pm\hbox{Id}$. 

The geodesics of $\mathbb H$ are half-circles (with respect to the usual
Euclidean metric of $\mathbb C$) centered at the boundary $\mathbb R$
of $\mathbb H\subset \mathbb C$
or halflines $\{a+iy\in \mathbb C\mid y>0\}\subset \mathbb H$ perpendicular to $\mathbb R$.

An orientation-preserving isometry $\iota$ of the Poincar\'e half-plane
is \emph{hyperbolic} if it admits an invariant geodesic on which it acts
by a translation. A M\"obius transformation 
associated to
$ \left(\begin{array}{cc}a&b\\c&d\end{array}\right)\in
\mathrm{SL}_2(\mathbb R)$ defines a hyperbolic isometry 
$\iota$ if and only if $\vert a+d\vert >2$.
The invariant geodesic of $\iota$
is given by the halfcircle in $\mathbb H$ delimited by the two 
real points
$\frac{a-d\pm\sqrt{(d-a)^2+4bc}}{2c}$ if $c\not=0$ respectively by
the halfline $\{b/(d-a)+iy\ \vert\ y>0\}$ otherwise. 

\begin{lem} \label{lemgeod}
For $\left(\begin{array}{cc}a&b\\c&d\end{array}\right)\in 
\mathrm{GL}_2(\mathbb R)$ with $cd\not=0$ and positive determinant 
$ad-bc>0$,
the image of the map from $\mathbb R_{>0}=(0,+\infty)$ 
into $\mathbb C$ defined by
$$t\longmapsto \frac{ait+b}{cit+d}$$
is an open half-circle of $\mathbb H$ centered on $\frac{ad+bc}{2cd}$
with radius $\left\vert\frac{ad-bc}{2cd}\right\vert$ (with respect to the 
Euclidean norm $\parallel z\parallel =\sqrt{x^2+y^2}$ for $z=x+iy\in 
\mathbb C$).
\end{lem}

\noindent{\bf Proof} M\"obius transformations 
preserve geodesics of $\mathbb H$.
Thus the M\"obius transformation defined by the
matrix $\left(\begin{array}{cc}a&b\\c&d
\end{array}\right)$ sends the geodesic $\{ti\ \vert\ 
t>0\}\subset\mathbb H$
onto a geodesic of $\mathbb H$ 
with finite boundary points given by $\frac{b}{d}$ (corresponding 
to $t=0$) and  $\frac{a}{c}$ (corresponding to $t=+\infty$).
This geodesic is the open halfcircle (for the usual Euclidean
metric of $\mathbb C$) of the upper halfplane
with center $\frac{1}{2}\left(\frac{b}{d}+\frac{a}{c}\right)=
\frac{ad+bc}{2cd}$ and diameter
$\left\vert\frac{a}{c}-\frac{ad+bc}{2cd}
\right\vert=\left\vert \frac{ad-bc}{2cd}\right\vert$.\hfill$\Box$

\begin{rem} Lemma \ref{lemgeod} is equivalent to the identity 
%$$\left\vert\frac{a+ibt}{c+idt}-\frac{ad+bc}{2cd}\right\vert^2=
$$\left(\frac{act^2+bd}{c^2t^2+d^2}-\frac{ad+bc}{2cd}\right)^2+\left(
\frac{ad-bc}{c^2t^2+d^2}\right)^2t^2=
\left(\frac{ad-bc}{2cd}\right)^2\ .$$
\end{rem}

%%%%%%%%%%%%%%%%%%%%%%%%%%%%%%%%%%%%%%%%%%%%%%%%%%%%%%%%%%%%%%%%%%%%%%%%%%

\section{Continued fractions}\label{sectcontfractions}

We denote by $[a_0;a_1,a_2,\dots]$ the continued fraction expansion
\begin{eqnarray*}
\theta&=&a_0+\frac{1}{a_1+\frac{1}{a_2+\ddots}}\ ,
\end{eqnarray*}
of a real number $\theta$. The coefficients $a_0,a_1,\dots$ are recursively 
defined by $a_i=\lfloor \theta_i\rfloor$
where $\theta_0=\theta$ and 
$\theta_n=\frac{1}{\theta_{n-1}-a_{n-1}}=\frac{1}{\lbrace \theta_{n-1}\rbrace}$ 
if $\theta_{n-1}\not\in \mathbb Z$, respectively by $\theta_n=0$
if $\theta_{n-1}\in \mathbb Z$. 
The coefficient $a_0$ of a continued fraction expansion can be an arbitrary
integer (positive, zero or negative). $a_1,a_2,\dots$
are either all strictly positive or they start with a finite number of
strictly positive integers followed by an infinite string of zeros. 
The last case arises if and only if $\theta$ is rational. 
The sequence $a_0,a_1,\dots$ is infinite and ultimately periodic 
with non-zero period
if and only if $\mathbb Q[\theta]$ is a quadratic number field.
Every irrational number has a unique continued fraction expansion.
Rational numbers have two expansions given by
$[a_0;a_1,\dots,a_m,1]$ and $[a_0;a_1,\dots,a_m+1]$
for suitable integers $m\geq 0,\ a_0\in \mathbb Z,\ a_1,\dots,a_m\geq 1$.

We have the continued fraction expansions
\begin{eqnarray}\label{defxn}
\theta_n&=&[a_n;a_{n+1},a_{n+2},a_{n+3},\dots]
\end{eqnarray}
for all $n\in\mathbb N$.
\emph{Convergents} for 
$$\theta=\theta_0=a_0+\frac{1}{a_1+\frac{1}{a_2+\ddots}}$$
are rational numbers of the form
$$\frac{p_{-2}}{q_{-2}}=\frac{0}{1},\ \frac{p_{-1}}{q_{-1}}=\frac{1}{0},
\frac{p_{n}}{q_{n}}=\frac{p_{n-2}+a_np_{n-1}}{q_{n-2}+a_nq_{n-1}}=
[a_0;a_1,\dots,a_n],\ n\geq 0$$
and can also be define by $\frac{p_n}{q_n}=[a_0;a_1,\dots,a_n]$, see
Theorem 149 in \cite{HW}.
\emph{Intermediate convergents} are given by
$$\frac{p_{n-2}+kp_{n-1}}{q_{n-2}+kq_{n-1}},\ k\in\{0,\dots,a_n-1\}\ .$$

The easy identity 
\begin{eqnarray}\label{pnqnidentity}
p_{n-1}q_{n}-p_nq_{n-1}=(-1)^n
\end{eqnarray}
(see Theorem 150 in \cite{HW} or Theorem 2 in \cite{Kh}), equivalent to 
$\frac{p_{n-1}}{q_{n-1}}-\frac{p_n}{q_n}=\frac{(-1)^n}{q_{n-1}q_n}$, implies
$$\frac{p_{2n}}{q_{2n}}<\frac{p_{2n+2}}{q_{2n+2}}<\dots<\theta<\dots<
\frac{p_{2n+1}}{q_{2n+1}}<\frac{p_{2n-1}}{q_{2n-1}}\ .$$
This shows
\begin{eqnarray}\label{contfracapprineq}
\left\vert\theta-\frac{p_n}{q_n}\right\vert< \frac{1}{a_{n+1}q_n^2}
\end{eqnarray} 
and ensures that convergents are excellent rational approximations
of an irrational number.

The following result is essentially identity 10.3.2 of \cite{HW}:

\begin{prop}\label{propx=reduite}
We have 
$$\theta=\frac{p_{n-2}+\theta_np_{n-1}}{q_{n-2}+\theta_nq_{n-1}}$$
for all $n\geq 0$.
\end{prop}

\noindent{\bf Proof} The result holds for $n=0$.

We have 
\begin{eqnarray*}
\frac{p_{n-1}+\theta_{n+1}p_n}{q_{n-1}+\theta_{n+1}q_n}&=&
\frac{p_{n-1}+\frac{1}{\theta_{n}-a_{n}}(p_{n-2}+a_{n}p_{n-1})}
{q_{n-1}+\frac{1}{\theta_{n}-a_{n}}(q_{n-2}+a_{n}q_{n-1})}\\
&=&\frac{p_{n-2}+\theta_{n}p_{n-1}}{q_{n-2}+\theta_{n}q_{n-1}}
\end{eqnarray*}
which ends the proof by induction.\hfill$\Box$ 

\begin{lem}\label{lemformthj} We have
\begin{eqnarray*}
\frac{\theta_{n+1}}{q_{n-1}+\theta_{n+1}q_n}&=&\frac{1}{q_{n-2}+\theta_nq_{n-1}}\ .
\end{eqnarray*}
\end{lem}

\noindent{\bf Proof}
We have
\begin{eqnarray*}
\frac{\theta_{n+1}}{q_{n-1}+\theta_{n+1}q_n}&=&\frac{1}{\frac{1}{\theta_{n+1}}q_{n-1}+q_n}\\
&=&\frac{1}{(\theta_n-a_n)q_{n-1}+q_n}\\
&=&\frac{1}{\theta_nq_{n-1}-a_nq_{n-1}+q_{n-2}+a_nq_{n-1}}
\end{eqnarray*}
where we have used the recursive definitions 
$\frac{1}{\theta_{n+1}}=\theta_n-a_n$ and
$q_n=q_{n-2}+a_nq_{n-1}$ of $\theta_{n+1}$ and of $q_n$.\hfill$\Box$

\begin{prop}\label{propthetapnx} We have
\begin{eqnarray}\label{propeqnid}
\theta-\frac{p_{n-2}+xp_{n-1}}{q_{n-2}+xq_{n-1}}=
\frac{(\theta_n-x)(-1)^n}{(q_{n-2}+xq_{n-1})(q_{n-2}+\theta_nq_{n-1})}\ .
\end{eqnarray}
\end{prop}

\noindent{\bf Proof} Proposition \ref{propx=reduite} shows that 
the result holds for $x=\theta_n$.

Since (\ref{propeqnid}) is equivalent to the identity
\begin{eqnarray}
(q_{n-2}+xq_{n-1})\theta-(p_{n-2}+xp_{n-1})
=\frac{\theta_n-x}{q_{n-2}+\theta_nq_{n-1}}(-1)^n
\end{eqnarray}
involving affine functions of $x$, it is enough to show the equality 
\begin{eqnarray}\label{affxeq}
q_{n-1}\theta-p_{n-1}=-\frac{(-1)^n}{q_{n-2}+\theta_nq_{n-1}}\ .
\end{eqnarray}
This holds for $n=0$ since it boils down to $-1=-1$. 
By induction, we have for $x=a_n$ the identity
$$\theta-\frac{p_{n-2}+a_np_{n-1}}{q_{n-2}+a_nq_{n-1}}=
\frac{(\theta_n-a_n)(-1)^n}{(q_{n-2}+a_nq_{n-1})(q_{n-2}+\theta_nq_{n-1})}$$
which can be rewritten as
\begin{eqnarray}\label{idthpnqn}
\theta-\frac{p_n}{q_n}=
\frac{(-1)^n}{\theta_{n+1}q_n(q_{n-2}+\theta_nq_{n-1})}
\end{eqnarray}
using the recursive definitions of $p_n,\ q_n$ and $\theta_{n+1}$.

The identity 
$$\theta_{n+1}\left(q_{n-2}+\theta_nq_{n-1}\right)=q_{n-1}+\theta_{n+1}q_n$$
equivalent to Lemma \ref{lemformthj} yields now (\ref{affxeq})
for $n+1$.\hfill$\Box$

\begin{rem} Identity (\ref{idthpnqn}) (corresponding to the specialization
$x=a_n$ of Proposition \ref{propthetapnx})
strengthens inequality (\ref{contfracapprineq})
since $a_{n+1}=\lfloor \theta_{n+1}\rfloor\leq \theta_{n+1}$
and $q_n=q_{n-2}+a_nq_{n-1}\leq q_{n-2}+\theta_nq_{n-1}$.
\end{rem}

\begin{lem}\label{lemcalcul} We have 
$$\frac{1}{q_{n-2}+\theta_{n}q_{n-1}}-\frac{a_{n+1}}{q_{n-1}+\theta_{n+1}q_n}=
\frac{1}{q_n+\theta_{n+2}q_{n+1}}\ .$$
\end{lem}

\noindent{\bf Proof} Using the identities
\begin{eqnarray*}
q_k&=&q_{k-2}+a_kq_{k-1},\\ 
\theta_{k+1}&=&\frac{1}{\theta_k-a_k},\\
\end{eqnarray*}
for $k=n$ and $k=n+1$ we have
\begin{eqnarray*}
&&\frac{1}{q_{n-2}+\theta_{n}q_{n-1}}-\frac{a_{n+1}}{q_{n-1}+\theta_{n+1}q_n}\\
&&=\frac{1}{q_n-a_{n}q_{n-1}+\theta_{n}q_{n-1}}-\frac{a_{n+1}}{q_{n-1}+\theta_{n+1}q_n }\\
&&=\frac{1}{q_n+\frac{1}{\theta_{n+1}}q_{n-1}}-\frac{a_{n+1}}{q_{n-1}+\theta_{n+1}q_n}\\
&&=\frac{\theta_{n+1}-a_{n+1}}{q_{n-1}+\theta_{n+1}q_n}\\
&&=\frac{\theta_{n+1}-a_{n+1}}{q_{n-1}+a_{n+1}q_n+(\theta_{n+1}-a_{n+1})q_n}\\
&&=\frac{1}{\theta_{n+2}q_{n+1}+q_n}
\end{eqnarray*}
which ends the proof.\hfill$\Box$

%%%%%%%%%%%%%%%%%%%%%%%%%%%%%%%%%%%%%%%%%%%%%%%%%%%%%%%%%%%%%%%%%%%%%%%%%%

\section{Linearization}\label{sectlinearization}

\begin{prop}\label{proplin} We have
%\begin{eqnarray}\label{formlineariz}
%\varphi_\theta(n+aq_{j-1}+bq_j)-\varphi_\theta(n)
%=\left(\frac{aq_{j-1}+bq_j}{2\sqrt{n}}+2i\pi\sqrt{n}\gamma\right)
%e^{2i\pi n\theta}+E
%\end{eqnarray}
\begin{eqnarray*}
\varphi_\theta(n+aq_{j-1}+bq_j)-\varphi_\theta(n)
&=&a\left(\frac{q_{j-1}}{2\sqrt{n}}-\frac{(-1)^j2i\pi\sqrt{n}\theta_{j+1}}
{q_{j-1}+\theta_{j+1}q_j}\right)
e^{2i\pi n\theta}\\
&&+b\left(\frac{q_{j}}{2\sqrt{n}}+\frac{(-1)^j2i\pi\sqrt{n}\theta_{j+2}}
{q_{j}+\theta_{j+2}q_{j+1}}\right)
e^{2i\pi n\theta}\\
&&+E_j(a,b)
\end{eqnarray*}
where $q_k$ is the denominator of the $k-$th
convergent $\frac{p_k}{q_k}=[a_0;a_1,\dots,a_k]$ for
$\theta=[a_0;a_1,a_2,\dots]$, where  
$\theta_k=[a_k;a_{k+1},a_{k+2},\dots]$ is defined as in (\ref{defxn})
and where the error $E_j(a,b)$ is asymptotically given by
\begin{eqnarray*}
&\left(\frac{-1}{8n^{3/2}}(aq_{j-1}+bq_j)^2+
\frac{i\pi}{\sqrt{n}}(aq_{j-1}+bq_j)\delta
-2\pi^2\sqrt{n}\delta^2\right)e^{2i\pi n\theta}
\end{eqnarray*}
with
$$\delta=-a\frac{(-1)^j\theta_{j+1}}{q_{j-1}+\theta_{j+1}q_j}+b
\frac{(-1)^j\theta_{j+2}}{q_j+\theta_{j+2}q_{j+1}}$$
if $\vert\varphi_\theta(n+aq_{j-1}+bq_j)-\varphi_\theta(n)\vert=O(1)$.
\end{prop} 

\begin{cor}\label{corlattice}
If $q_{j-1}$ and $q_j$ are denominators of two consecutive 
convergents 
$\frac{p_{j-1}}{q_{j-1}}$ and $\frac{p_j}{q_j}$ of $\theta$ 
such that $q_{j-1}\leq \sqrt n<q_j$, then the smallest points of 
$$\left(\varphi_\theta(\mathbb N)-\varphi_\theta(n)\right)e^{-2i\pi n\theta}$$
are close to the smallest points of the lattice 
\begin{eqnarray}\label{linlattice}
\mathbb Z\left(\frac{q_{j-1}}{2\sqrt{n}}-(-1)^j\frac{2i\pi\sqrt{n}\theta_{j+1}}{q_{j-1}+\theta_{j+1}q_j}\right)+
\mathbb Z\left(\frac{q_j}{2\sqrt{n}}+(-1)^j\frac{2i\pi\sqrt{n}\theta_{j+2}}{q_j+\theta_{j+2}q_{j+1}}\right)
\end{eqnarray}
with an error of order $O\left(\frac{1}{\sqrt{n}}\right)$.
\end{cor}

The lattice described by (\ref{linlattice}) contains always a non-zero
element of absolute value smaller than $\sqrt{\frac{1}{4}+4\pi^2}<3\pi$.

\begin{rem} Fundamental domains of the lattice $\Lambda$ given by 
(\ref{linlattice}) have area $\pi$ as shown by the 
identities
\begin{eqnarray*}&&(-1)^j\det\left(\begin{array}{cc}
\frac{q_{j-1}}{2\sqrt{n}}&-\frac{(-1)^j2i\pi\sqrt{n}\theta_{j+1}}{q_{j-1}+
\theta_{j+1}q_j}\\
\frac{q_{j}}{2\sqrt{n}}&+\frac{(-1)^j2i\pi\sqrt{n}\theta_{j+2}}{q_j+
\theta_{j+2}q_{j+1}}\end{array}\right)\\
&=&\left(q_{j-1}\frac{\theta_{j+2}}{q_j+\theta_{j+2}q_{j+1}}+
q_j\frac{\theta_{j+1}}{q_{j-1}+\theta_{j+1}q_j}\right)\pi\\
&=&\left(q_{j-1}\frac{1}{q_{j-1}+\theta_{j+1}q_j}+
q_j\frac{\theta_{j+1}}{q_{j-1}+\theta_{j+1}q_j}\right)\pi\\
&=&\pi
\end{eqnarray*}
where the second equality is given by Lemma \ref{lemformthj}.

Since the regular hexagonal lattice has maximal density,
the lattice $\Lambda$ contains always a non-zero element of absolute 
value at most $\sqrt{\frac{2\pi}{\sqrt{3}}}\sim 1.9046$.
\end{rem}

\noindent{\bf Proof of Proposition \ref{proplin}} Setting
$$F(s,\gamma)=\sqrt{n+s}\ e^{2i\pi n\theta+2i\pi\gamma}$$
we want to approximate
$$F(aq_{j-1}+bq_j,\delta)-F(0,0)$$
where $\delta=(aq_{j-1}+bq_j)\theta-c$ is the difference between 
$(aq_{j-1}+bq_j)\theta$ and the integer $c$ closest to $(aq_{j-1}+bq_j)\theta$.
We do this in the usual way by considering the linear approximation 
$$L=\frac{\partial F}{\partial s}(0,0)(aq_{j-1}+bq_j)+
\frac{\partial F}{\partial \gamma}(0,0)\delta$$
and by estimating the error using second-order derivatives.
The necessary partial derivatives of $F$ are:
\begin{eqnarray*}
%F(0,0)&=&\sqrt{n}e^{2i\pi n\theta},\\
\frac{\partial F}{\partial s}(0,0)&=&\frac{1}{2\sqrt{n}}e^{2i\pi n\theta},\\
\frac{\partial F}{\partial \gamma}(0,0)&=&2i\pi\sqrt{n}e^{2i\pi n\theta},\\
\frac{\partial^2 F}{\partial s^2}(0,0)&=&\frac{-1}{4n^{3/2}}e^{2i\pi n\theta},\\
\frac{\partial^2 F}{\partial s\partial \gamma}(0,0)&=&
\frac{i\pi}{\sqrt{n}}e^{2i\pi n\theta},\\
\frac{\partial^2 F}{\partial \gamma^2}(0,0)&=&-4\pi^2\sqrt{n}e^{2i\pi n\theta}.\\
\end{eqnarray*}
The contribution coming from $\frac{\partial F}{\partial s}(0,0)(aq_{j-1}+bq_j)$
to $L$ is given by 
$$\frac{1}{2\sqrt{n}}e^{2i\pi n\theta}(aq_{j-1}+bq_j)\ .$$
In order to compute $\frac{\partial F}{\partial \gamma}(0,0)\delta$ 
we split $\delta$ into $\delta=a\delta_{j-1}+b\delta_j$
where $\delta_k$ for $k\in\{j-1,j\}$
is the difference between $q_k\theta$ and the integer 
closest to $q_k\theta$. Since $q_k$ is a denominator of the convergent 
$\frac{p_k}{q_k}=[a_0;a_1,\dots,a_k]$ of $\theta$, this integer is given by the numerator $p_k$. We have
\begin{eqnarray*}
\theta-\frac{p_k}{q_k}&=&\frac{p_k+\theta_{k+2}p_{k+1}}{q_k+\theta_{k+2}q_{k+1}}
-\frac{p_k}{q_k}\\
&=&\frac{(p_{k+1}q_k-p_kq_{k+1})\theta_{k+2}}{q_k(q_k+\theta_{k+2}q_{k+1}}\\
&=&(-1)^k\frac{\theta_{k+2}}{q_k(q_k+\theta_{k+2}q_{k+1})}
\end{eqnarray*}
where we have used Proposition \ref{propx=reduite} and 
identity (\ref{pnqnidentity}). This yields
$$\delta_k=(-1)^k\frac{\theta_{k+2}}{q_k+\theta_{k+2}q_{k+1}}$$
and shows
\begin{eqnarray*}
\frac{\partial F}{\partial s}(0,0)\delta&=&
\frac{\partial F}{\partial s}(0,0)(a\delta_{j-1}+b\delta_j)\\
&=&\frac{\partial F}{\partial s}(0,0)
\left(-a\frac{(-1)^j\theta_{j+1}}{q_{j-1}+\theta_{j+1}q_{j}}+
b\frac{(-1)^j\theta_{j+2}}{q_{j}+\theta_{j+2}q_{j+1}}\right).
\end{eqnarray*}
The order of the error is given by
$$\frac{1}{2}\frac{\partial^2 F}{\partial s^2}(0,0)(aq_{j-1}+bq_j)^2+
\frac{\partial^2 F}{\partial s\partial\gamma}(0,0)(aq_{j-1}+bq_j)\delta+
\frac{1}{2}\frac{\partial^2 F}{\partial \gamma^2}(0,0)\delta^2$$
and can be evaluated easily.\hfill$\Box$

%%%%%%%%%%%%%%%%%%%%%%%%%%%%%%%%%%%%%%%%%%%%%%%%%%%%%%%%%%%%%%%%%%%%%%%%%%%

\section{The phyllotactic geodesic and proof of Theorem \ref{mainthm}}\label{sectphyllgeod}

Using (\ref{lemformthj}) we can rewrite the lattice $\Lambda$ given by formula
(\ref{linlattice}) of Corollary \ref{corlattice} as
$$\mathbb Z\left(\frac{q_{j-1}}{2\sqrt{n}}-(-1)^j\frac{2i\pi\sqrt{n}}{q_{j-2}+\theta_jq_{j-1}}\right)+
\mathbb Z\left(\frac{q_j}{2\sqrt{n}}+(-1)^j\frac{2i\pi\sqrt{n}}{q_{j-1}+\theta_{j+1}q_{j}}\right)\ .$$
In particular, the lattice $\Lambda$ is similar to 
the lattice $\mathbb Z+\mathbb Z \tau_j(n)$ where 
\begin{eqnarray}\label{deftau}
\tau_j(t)=-(-1)^j
\frac{q_{j-1}-(-1)^j\frac{4i\pi t}{q_{j-2}+\theta_jq_{j-1}}}
{q_{j}+(-1)^j\frac{4i\pi t}{q_{j-1}+\theta_{j+1}q_j}}\ .
\end{eqnarray}

\begin{rem} A straightforward computation shows that the imaginary 
part of $\tau_j(t)$, given by 
$$\frac{\frac{q_{j-1}}{q_{j-1}+\theta_{j+1}q_j}+\frac{q_j}{q_{j-2}+\theta_j
q_{j-1}}}{q_j^2+\frac{16\pi^2t^2}{\left(q_{j-1}+\theta_{j+1}q_j\right)^2}}
4\pi t\ ,$$
is strictly positive if $t$ is strictly positive.
\end{rem}

\begin{thm}\label{thmtaujeq}
For all $j\geq 0$ we have 
$$(\tau_j(t)-(-1)^ja_{j+1})\tau_{j+1}(t)=-1\ .$$
\end{thm}

\noindent{\bf Proof} Theorem \ref{thmtaujeq} boils down to the identity
\begin{eqnarray*}
&&q_{j-1}-(-1)^j\frac{4i\pi t}{q_{j-2}+\theta_jq_{j-1}}+a_{j+1}q_j+a_{j+1}(-1)^j
\frac{4i\pi t}{q_{j-1}+\theta_{j+1}q_j}\\
&=&q_{j+1}-(-1)^j
\frac{4i\pi t}{q_j+\theta_{j+2}q_{j+1}}\ .
\end{eqnarray*}
The identity $q_{j+1}=q_{j-1}+a_{j+1}q_j$ shows that 
the constant parts (with respect to $t$) 
of both sides are equal.
Linear coefficients of $t$ are equal by Lemma \ref{lemcalcul}.\hfill$\Box$

\noindent{\bf Proof of Theorem \ref{mainthm}}
By Theorem \ref{thmtaujeq}, the two geodesics defined by
$\tau_j$ and $\tau_{j+1}$ are related by the integral M\"obius
transformations
\begin{eqnarray}\label{homAjAjp}
\tau_{j+1}=\frac{-1}{\tau_j-(-1)^ja_{j+1}}=\left(\begin{array}{cc}
0&1\\-1&(-1)^ja_{j+1}\end{array}\right) \cdot \tau_j
\end{eqnarray}
and
\begin{eqnarray}\label{homAjpAj}
\tau_j=\frac{(-1)^ja_{j+1}\tau_{j+1}-1}{\tau_{j+1}}=\left(\begin{array}{cc}
(-1)^ja_{j+1}&-1\\1&0\end{array}\right) \cdot \tau_{j+1}\ .
\end{eqnarray}
Thus they project onto a unique geodesic 
on the modular curve $\mathrm{PSL}_2(\mathbb Z)\backslash\mathbb H$
represented for example by
\begin{eqnarray}\label{formphyllgeod}
\tau_0(t)=-\frac{q_{-1}-\frac{4i\pi t}{q_{-2}+\theta_0q_{-1}}}{q_0+
\frac{4i\pi t}{q_{-1}+\theta_1q_0}}=
\frac{4i\pi t}{1+\frac{4i\pi t}{\theta_1}}=
\frac{4i\pi t}{1+4i\pi \lbrace\theta\rbrace t}
\end{eqnarray}
(where $q_{-2}=1,\ q_{-1}=0,\ q_0=1$ and 
$\theta_1=\frac{1}{\theta-a_0}=\frac{1}{\lbrace \theta\rbrace}$).
 
This implies Theorem \ref{mainthm} since the linearization error 
is of order $O\left(\frac{1}{\sqrt{n}}\right)$
for elements of $\varphi_\theta(\mathbb N)$ at bounded distance from 
$\varphi_\theta(n)$.\hfill$\Box$

Formula (\ref{formphyllgeod}) defines a geodesic of the hyperbolic 
half-plane for every real number $\theta$. Indeed, (\ref{formphyllgeod})
is a vertical half-line (and thus a geodesic) if $\theta$ is integral
and it defines a halfcircle of $\mathbb H$ orthogonal to $\mathbb R$
(and thus a geodesic) with boundary points $0$ corresponding to 
$t=0$ and $\frac{1}{\{\theta\}}$ corresponding to $t=\infty$
otherwise.

Diophantine properties of $\theta$ are related to the dynamical behaviour
of the geodesic $\gamma_\theta$ projected onto $
\mathrm{PSL}_2(\mathbb Z)\backslash\mathbb H$ as follows: 
after starting at the cusp,
(the projection of) 
$\gamma_\theta$ turns (slightly less) than $a_1$ times around the cusp
before passing between the two conical points of the modular curve. It turns
then in the same sense (and slightly less than) $a_2$ times around the cusp before crossing
again the shortest geodesic segment joining the two conical points and so on.
A large coefficient $a_k$ causes the (projection of the) geodesic 
$\gamma_\theta$ to climb the modular curve up to a height given asymptotically
(in $a_k$) by $a_ki\in\mathcal M$. 
This gives rise to points of $\varphi_\theta(\mathbb N)$
having Voronoi cells with bad isoperimetric properties.

For a divergence angle 
$2\pi\theta$ determined by the golden ratio $\theta=\frac{1+\sqrt{5}}{2}$
(or close relatives of it) the continued fraction expansion involves 
only ones (or only ones after perhaps a few initial ``accidents'').
This is the optimal situation leaving no possibility of improvement.
In particular, 
the phyllotactic geodesic $\gamma_{(1+\sqrt{5})/2}$ is asymptotically
equal to the geodesic 
$$t\longmapsto \tilde \gamma(t)=\frac{(-1+\sqrt{5})it-1-\sqrt{5}}{2(it+1)}
=\frac{-1-\sqrt{5}+(-1+\sqrt{5})t^2+2it\sqrt{5}}{2(1+t^2)}$$ 
with boundary points $\frac{-1\pm\sqrt{5}}{2}$
(and containing the points $-1+i=\tilde\gamma\left(\frac{-1+\sqrt{5}}{2}
\right)$ and $i=\tilde\gamma\left(\frac{1+\sqrt{5}}{2}\right)$ of $\mathbb H$).
The equality
$$\tilde \gamma(\theta^4t)=\frac{\tilde\gamma(t)+1}
{\tilde\gamma(t)+2}$$ shows that $\tilde\gamma$ is invariant under
the integral M\"obius transformation defined by the matrix 
$\left(\begin{array}{cc}1&1\\1&2\end{array}\right)=
\left(\begin{array}{cc}0&1\\1&1\end{array}\right)^2$. It projects onto
the shortest closed geodesic of the modular curve $
\mathrm{PSL}_2(\mathbb Z)\backslash\mathbb H$. 

Since $\frac{d}{dt}\tilde\gamma(t)=\frac{-\sqrt{5}}{(t-i)^2}$, the 
parametrized geodesic $\tilde\gamma(t)$ has the same instant speed
$\frac{1}{\Im(\tilde\gamma(t))}\vert \tilde\gamma'(t)\vert=\frac{1}{t}$
as the phyllotactic geodesic $\gamma_\theta(t)$.
We have thus asymptotically $\gamma_\theta(t)\sim\gamma_\theta\left
(\theta^4t\right)$. The phyllotactic set
$\varphi_{\theta}(\mathbb N)$ has thus almost isometrical neighbourhoods
around $\varphi_\theta(n)$ and $\varphi_\theta(m)$ if $n$ is large and $m$ is 
close to $\theta^4n$.

%%%%%%%%%%%%%%%%%%%%%%%%%%%%%%%%%%%%%%%%%%%%%%%%%%%%%%%%%%%%%%%%%%%%%%%%%%

\section{A geometric construction}\label{sectgeomconstr}

We construct in this section a slightly 
different geodesic on the modular domain
$\mathrm{PSL}_2(\mathbb Z)\backslash\mathbb H$ which is asymptotically
associated to linearizations of  $\varphi_\theta(\mathbb N)$.

We denote by $\mathcal L=\mathcal L_\theta$ the line $\mathbb R(1,-\theta)$
of slope $-\theta$ containing the origin.
A convergent $\frac{p_j}{q_j}=[a_0;a_1,\dots,a_j]$ of 
$\theta=[a_0;a_1,a_2,\dots]$ yields an integral
point $(q_j,-p_j)$ close to $\mathcal L$ as follows: Denoting by 
$\pi_{\mathcal L}(q_j,-p_j)$ the orthogonal projection
of $(q_j,-p_j)$ onto $\mathcal L$ we have
$$\pi_{\mathcal L}(q_j,-p_j)=q_j(1,-\theta)+O\left(\frac{1}{q_j}\right)
(1,-\theta)$$
for all $j\geq 0$. The equality $q_j\theta-p_j=
(-1)^j\frac{\theta_{j+2}}{q_j+\theta_{j+2}q_{j+1}}$
corresponding to the case  $n=j+2$ and $x=0$ of 
Proposition \ref{propthetapnx} implies the identity 
$$(q_j,-p_j)-q_j(1,-\theta)=
\left(0,(-1)^j\frac{\theta_{j+2}}{q_j+\theta_{j+2}q_{j+1}}\right)\ .$$
The orthogonal projection $\pi_{\mathcal L^\perp}(q_j,-p_j)$ 
of $(q_j,-p_j)$ onto the line $\mathcal L^{\perp}=\mathbb R(\theta,1)$ orthogonal
to $\mathcal L=\mathbb R(1,-\theta)$ is given by
\begin{eqnarray*}
\pi_{\mathcal L^\perp}(q_j,-p_j)&=&\pi_{\mathcal L^\perp}\big((q_j,-p_j)-q_j(1,-\theta)
\big)\\
&=&\frac{(-1)^j\theta_{j+2}}{(q_j+\theta_{j+2}q_{j+1})}\frac{(\theta,1)}{(1+\theta^2)}\ .
\end{eqnarray*}
Thus we can rewrite the right side of the obvious identity
\begin{eqnarray*}
(q_j,-p_j)&=&\pi_{\mathcal L}(q_j,-p_j)+\pi_{\mathcal L^\perp}(q_j,-p_j)
\end{eqnarray*}
as 
\begin{eqnarray}\label{projqjpmj}
q_j(1,-\theta)+
\frac{(-1)^j4\pi n\theta_{j+2}}{(q_j+\theta_{j+2}q_{j+1})}\frac{(\theta,1)}
{4\pi n(1+\theta^2)}+O\left(\frac{1}{q_j}\right)(1,-\theta)\ .
\end{eqnarray}

We endow now $\mathbb R^2$ with an Euclidean metric $ds_{\theta,n}$
turning the vectors $(1,-\theta),\frac{(\theta,1)}{4\pi n(1+\theta^2)}$ 
into an orthogonal basis. Comparision of (\ref{projqjpmj}) with
(\ref{linlattice}) shows that the lattice 
$(\mathbb Z^2,ds_{\theta,n})$ is asymptotically equivalent with the 
linearization of $\varphi_\theta(\mathbb N)$ at the point 
$\varphi_\theta(n)$. Since
\begin{eqnarray*}
(1,0)&=&\frac{1}{1+\theta^2}(1,-\theta)+4\pi n \theta
\frac{(\theta,1)}{4\pi n(1+\theta^2)},\\
(0,1)&=&-\frac{\theta}{1+\theta^2}(1,-\theta)+4\pi n
\frac{(\theta,1)}{4\pi n(1+\theta^2)},\\
\end{eqnarray*}
the Euclidean lattice $(\mathbb Z^2,ds_{\theta,n})$
corresponds to the point of $\mathrm{PSL}_{2}(\mathbb Z)\backslash\mathbb H$
represented by
\begin{eqnarray*}
&&\frac{1+4i\pi n\theta(1+\theta^2)}
{\theta-4i\pi n(1+\theta^2)}\\
&=&\frac{\theta-16\pi^2n^2\theta(1+\theta^2)^2+4i\pi
n(1+\theta^2)^2}{\theta^2+16\pi^2n^2(1+\theta^2)^2}\ .
\end{eqnarray*}
All these points are elements of the hyperbolic geodesic
with boundary points $\frac{1}{\theta}$ (for $n=0$) and $-\theta=
\left(\begin{array}{cc}a_0&1\\-1&0\end{array}\right)\frac{1}{
\lbrace \theta\rbrace}$ (for $n=\infty$).

\begin{rem} Since we have
\begin{eqnarray*}
&&\frac{1+4i\pi t\theta(1+\theta^2)}
{\theta-4i\pi t(1+\theta^2)}-
\left(\begin{array}{cc}a_0&1\\-1&0\end{array}\right)
\frac{4i\pi t}{4i\pi\lbrace \theta\rbrace t+1}\\
&=&\frac{1+4i\pi t\theta(1+\theta^2)}
{\theta-4i\pi t(1+\theta^2)}+\frac{4i\pi \theta t+1}{4i\pi t}\\
&=&\frac{\theta}{4i\pi t(\theta-4i\pi t(1+\theta^2))}\\
&=&\frac{1+\theta^2}{\theta^2+16\pi^2t^2(1+\theta^2)^2}-
\frac{\theta^2}{4\pi  t\left(\theta^2+16\pi^2t^2(1+\theta^2)^2\right)}i\ ,
\end{eqnarray*}
the hyperbolic distance between the two points 
$$\frac{1+4i\pi t\theta(1+\theta^2)}
{\theta-4i\pi t(1+\theta^2)}\hbox{ and }
\left(\begin{array}{cc}a_0&1\\-1&0\end{array}\right)
\frac{4i\pi t}{4i\pi\lbrace \theta\rbrace t+1}$$
of $\mathbb H$ is roughly given by
\begin{eqnarray*}
&&\frac{\theta^2+16\pi^2t^2(1+\theta^2)^2}{4\pi t(1+\theta^2)^2}
\left\vert\frac{1+\theta^2}{\theta^2+16\pi^2t^2(1+\theta^2)^2}-
\frac{\theta^2}{4\pi  t\left(\theta^2+16\pi^2t^2(1+\theta^2)^2\right)}i
\right\vert
\end{eqnarray*}
which simplifies to
\begin{eqnarray*}
\left\vert\frac{1}{4\pi t(1+\theta^2)}-\frac{\theta^2}{16\pi^2t^2}i
\right\vert\ .
\end{eqnarray*}
Thus it is asymptotically equal to 
$\frac{1}{4\pi t(1+\theta^2)}$ which is asymptotically 
much smaller than the error
$O\left(\frac{1}{\sqrt{n}}\right)$ due to linearization at a point
$\varphi_\theta(n)=\sqrt{n}e^{2i\pi \theta n}$ of order $O(\sqrt{t})$.
\end{rem}

%%%%%%%%%%%%%%%%%%%%%%%%%%%%%%%%%%%%%%%%%%%%%%%%%%%%%%%%%%%%%%%%%%%%%%%%%%%

\section{Metric properties of phyllotactic sets}\label{sectmetrprop}

A subset $S$ of a metric space $E$ is {\em uniformly discrete} 
if there exists
a strictly positive real constant $\delta$ such that 
$d(a,b)\geq \delta$ for every pair $a,b$ of distinct points
in $S$. Equivalently, $S$ is uniformly discrete if 
open balls of radius $\delta/2$ centered at all elements of $E$ 
are disjoint (for a small strictly positive constant $\delta$).

A subset $\mathcal S$ of a metric space $E$ is an \emph{$\epsilon-$net}
if every point of $E$ is at distance at most $\epsilon$ from a point 
of $\mathcal S$. Equivalently, $E$ is covered by 
the set of closed balls of radius $\epsilon$ centered at 
elements of $\mathcal S$.

The following result is a straightforward consequence of the fact that 
are{\ae} of Voronoi cells defined by $\varphi_\theta(\mathbb N)$
are asymptotically equal to 
$\pi$ if $\theta$ is irrational:

\begin{prop} \label{propunifdens} The following assertions are equivalent:

(i) $\theta$ is irrational and has bounded coefficients $a_1,a_2,\dots$
in its continued fraction expansion $\theta=[a_0;a_1,a_2,\dots]$.

(ii) $\varphi_\theta(\mathbb N)$ is uniformly discrete in $\mathbb C$
(identified with the Euclidean plane in the obvious way).

(iii) $\varphi_\theta(\mathbb N)$ is an $\epsilon-$net of $\mathbb C$.

(iv) All Voronoi cells of $\varphi_\theta(\mathbb N)$ have bounded diameter.

(v) Discs of radius $R$ (and arbitrary centers) 
in $\mathbb C$ contain $R(R+O(1))$ 
points of $\varphi_\theta(\mathbb N)$.

(vi) The image $\gamma_\theta([1,\infty))$ of the phyllotactic \lq\lq 
half-geodesic'' is contained in a compact subset of the modular curve
$\mathrm{PSL}_2(\mathbb Z)\backslash\mathbb H$.
\end{prop}

We leave the proof to the reader.\hfill$\Box$

Observe that $\varphi_\theta(\mathbb N)$ is never uniformly discrete
if $\theta$ is rational.

%%%%%%%%%%%%%%%%%%%%%%%%%%%%%%%%%%%%%%%%%%%%%%%%%%%%%%%%%%%%%%%%%%%%%%%%%%%%

\section{Parastichy spirals}\label{sectparast}

We denote by $\partial \mathcal M$ the boundary in $\mathbb C$ 
of the fundamental domain $\mathcal M$ defined by (\ref{deffunddomM}). 
The interior
$\mathcal M\setminus \partial \mathcal M$ corresponds to lattices having 
a unique pair $\pm u$ of opposite shortest non-zero vectors and a unique
pair $\pm v$ of shortest vectors which are $\mathbb R-$linearly 
independent from $\pm u$. More precisely, for 
$z\in\mathcal M\setminus\partial{\mathcal M}$ the unique pair $\pm u$
of non-zero shortest vectors in $\mathbb Z+\mathbb Zz$ is given by $\pm 1$
and the unique pair $\pm v$ of shortest vectors outside
$\mathbb R$ coincides with $\pm z$. Notice that $\mathbb R-$linear
independency of $v$ from $u$ is necessary in order to discard
$\pm 2u,\pm 3u,\dots$ which might be smaller than $v$ for lattices
associated to $z\in\mathcal M$ with large modulus. 

Lattices corresponding to elements $z$ of norm $1$ in $\mathcal M$
have (at least) two pairs of shortest vectors given by $\pm 1$ and 
$\pm z$ in $\mathbb Z+\mathbb Zz$. The regular hexagonal 
lattice corresponding to $z=
\frac{1+i\sqrt{3}}{2}$ is the unique lattice with three 
pairs $\pm 1,\pm \frac{1+i\sqrt{3}}{2},\pm \frac{1-i\sqrt{3}}{2}$
of shortest non-zero 
vectors. Lattices associated to 
$z=\frac{1+it}{2}$ for $t>\sqrt{3}$ have a unique
pair $\pm 1$ of shortest non-zero-vectors and two non-real pairs 
$\pm z$ and $\pm (z-1)$ of shortest non-real vectors.

Connecting points of a lattice $\Lambda$ 
indexed by $z\in\mathcal M\setminus\partial{
\mathcal M}$ with their closest neighbours we get a
set of parallel lines. Joining closest lattice-points on two
such adjacent lines we obtain a second set of parallel lines.
These two sets of parallel lines cut the complex plane into 
fundamental domains for $\Lambda$ given by isometric rhombi.

\emph{Parastichy spirals} are analogues of these lines 
in $\varphi_\theta(\mathbb N)$. 
More precisely, we define (generically) the
\emph{primary parastichy spirals} of $\varphi_\theta(\mathbb N)$
as the piecewise-wise linear paths obtained by joining vertices of 
$\varphi_\theta(\mathbb N)$ to their two approximatively
opposite nearest neighbours. Similarly, we construct \emph{secondary 
parastichy spirals} by joining vertices of $\varphi_\theta(\mathbb N)$ 
to their nearest neighbours on adjacent neighbouring primary parastichy
spirals. 

Primary parastichy spirals exist essentially at every 
point far from the origin 
except where they become blurred with secondary parastichy spirals.
At such points (corresponding to crossings of the phyllotactic 
geodesic with the image of the unit circle in 
$\mathrm{PSL}_2(\mathbb Z)\backslash\mathbb H$),
primary and secondary parastichy spirals get exchanged.
We call such a situation a \emph{parastichy transition of type I}. 

Secondary parastichy spirals are however well-defined only
if the local situation corresponds to a lattice indexed by an element
of $\mathcal M$ which is not too close to the cusp. 
For example, the phyllotactic set $\varphi_\theta(\mathbb N)$ associated to 
a rational number $\theta=\frac{p}{q}$
is contained in $q$ half-rays originating at $0$. 
Far from the origin, primary parastichies (and are no longer spirals) 
coincide with 
these half-rays and secondary parastichies
make no longer sense. Moreover, for points
$z\in\mathcal M$ with real part close to $1/2$ (or $-1/2$) 
a family of secondary parastichy spirals fades
away and is replaced by a new family of secondary parastichy spirals,
giving rise to a \emph{parastichy transition of type II}.
A coefficient $a_i>1$ yields $a_i-1$ parastichy transitions of 
type II. The occurence of parastichy transitions of type II is easy to detect
visually: it leads
to much less uniform point distributions in $\varphi_\theta(\mathbb N)$.
Figure 4 displays two examples.

The geometric construction of Section \ref{sectgeomconstr} shows the 
well-known fact that primary and secondary parastichies
form two sets of spirals with different 
orientations if $\theta$ is irrational. 
Indeed, primary, respectively secondary, 
parastichies around $\varphi_\theta(n)$ are
defined by $\varphi_\theta(n\pm q_j)$, respectively by $\varphi_\theta(n
\pm(q_{j-1}+kq_j))$ for suitable integers $j,k$, see Figure 2
where one has to think of $(0,0)$ as the point $\varphi_\theta(n)$ 
and of $\mathcal L$ as the ray defined by $\mathbb R_{>0}\varphi_\theta(n)$.
(As always, $q_{j-1}$ and $q_j$ are denominators of convergents for $\theta$.) 
The same integers $j,k$ work for all $n$ in some interval
of large length compared to $\sqrt{n}$. 
The plane $\mathbb R^2$ is increasingly
squezed (for increasing $n$) in the direction of $\mathcal L$
and expanded in the orthogonal direction $\mathcal L^\perp$ in the 
construction of Section \ref{sectgeomconstr}. 
This implies that parastichies of both kinds 
bend away from the rays issued by the origin.
Moreover, there is exactly one parastichy family of larger,
respectively of smaller slope than $\varphi_\theta(n)$ as can 
be seen by inspecting Figure 2. This explains the apparition of 
crisscrossing spirals in Figure 1. Secondary parastichy families are 
however no longer discernible (to the eye) for larger values of $k$,
see Figure 4 where there are regions without obvious secondary parastichies.
For irrational $\theta$, they can however always be drawn ``by continuity'',
if we start in suitable regions where no problems occur and if we
push them forward using type II transition for (the projection of)
$\gamma_\theta$ crossing 
the infinite boundary segment of $\mathcal M$.

\subsection{Transitions for parastichy families}

The number of \lq\lq parallel'' 
primary parastichy spirals forming a common family is always
a denominator $q_j$ of a convergent $\frac{p_j}{q_j}=
[a_0;a_1,\dots,a_j]$ of the divergence angle $2\pi\theta=
2\pi[a_0,a_1,a_2,\dots]$.
The number of secondary parastichy spirals in a common family is a 
denominator $q_j+kq_{j+1},k\in \{0,1,\dots,a_{j+2}\}$
of an intermediate convergent.  

The exterior region of Figure 1 for example 
contains 55 primary parastichy spirals turning clockwise
and 34 secondary parastichy spirals turning counterclockwise.

The evolution of the numbers of parastichy spirals (PS in the following 
table) can be described
by:
$$\begin{array}{|c|c|c|c|}
&\vdots&\vdots&\vdots\\
\hline
&q_{j-1}\hbox{ prim. PS}&&q_j\hbox{ sec. PS}\\
I&&\vert \overline{\gamma_\theta(n)}\vert\sim 1&\\
&q_{j-1}\hbox{ sec. PS}&&q_j\hbox{ prim. PS}\\
\hline
II&q_{j-1}+q_j\hbox{ sec. PS}&\vert\Re(\overline{\gamma_\theta(n)}
\vert\sim \frac{1}{2}&q_j\hbox{ prim. PS}\\
\hline
II&q_{j-1}+2q_j\hbox{ sec. PS}&\vert\Re(\overline{\gamma_\theta(n)}
\vert\sim \frac{1}{2}&q_j\hbox{ prim. PS}\\
\hline
II&\vdots&\vdots&\vdots\\
\hline
&q_{j-1}+a_{j+1}q_j=q_{j+1}\hbox{ sec. PS}&&q_j\hbox{ prim. PS}\\
I&&\vert \overline{\gamma_\theta(n)}\vert\sim 1&\\
&q_{j+1}\hbox{ prim. PS}&&q_j\hbox{ sec. PS}\\
\hline
II&q_{j+1}\hbox{ prim. PS}&\vert\Re(\overline{\gamma_\theta(n)}
\vert\sim \frac{1}{2}&q_j+q_{j+1}\hbox{ sec. PS}\\
\hline
&\vdots&\vdots&\vdots
\end{array}$$
where $\overline{\gamma_\theta(n)}$ denotes a representant of 
of $\gamma_\theta(n)$ in $\mathcal M$.

\begin{figure}[h]\label{figure2}
\input{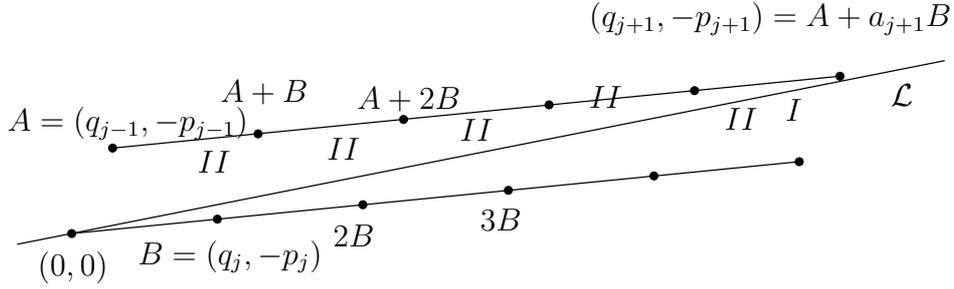}
\caption{Parastichy transitions}
\end{figure}

(Families of) parastichy spirals exist in some sense in the real world
(eg. in the approximate point set given by seeds of a real sunflower).
Parastichy transitions however are an ideal (or platonic) concept.
They  are in some sense ``smeared out'' (like the boundary between
adjacent colours of a rainbow) and cannot be localized exactly in a 
real flower.

Figure 2 attempts to illustrate 
the occurence of parastichy transitions using the geometric construction.

\begin{figure}[h]\label{figure3}
\epsfysize=6cm
\center{\epsfbox{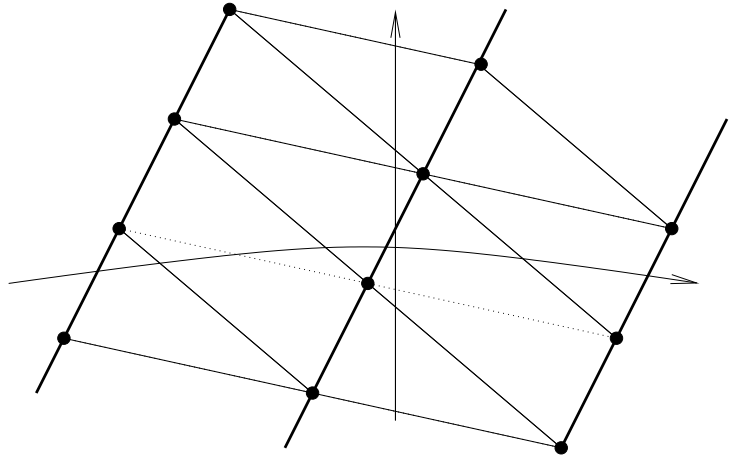}}
\caption{Death and birth of a family of secondary parastichy spirals}
\end{figure}
Figure 3 shows the death and birth of a family of secondary parastichy 
spirals corresponding to a value of $n$ such that 
$\vert\Re(\overline{\gamma_\theta(n)}\vert\sim \frac{1}{2}$.
The vertical arrow represents a ray issued from the origin.
Primary parastichies are represented by fat segments, the dying family 
of secondary parastichies is drawn with dotted segments
and the newborn family of secondary parastichies is given by 
ordinary segments. Turning around the origin on a circle of constant 
radius $\sim \sqrt{n}$ and counting the numbers $a$ of primary 
parastichy spirals, $b$ of dying secondary para\-stichies and $c$ 
of newborn secondary parastichies, one gets the relation $c=a+b$.
More precisely, the piecewise linear path involving only segments of 
primary and secondary dying para\-stichies giving the best approximation
of the circle with radius $\sqrt{n}$ consists of $b$
segments on primary parastichies and $a$ segments on secondary parastichies.
In order to work with the family of newborn parastichies, one has 
to replace every segment of the dying family by two segments, one from
a primary parastichy and one from a newborn secondary parastichy.
The number of segments on secondary parastichies (which is equal to 
the number of curves in the primary family) remains thus constant (and 
equals $a$) and the number $c$ of segments on primary parastichies
(which equals the number $c$ of curves in the newborn secondary
parastichy family) increases by $a$ to $c=a+b$.

In order to prove that parastichy families are enumerated by 
denominators of (intermediate) convergents, it is now enough to 
remark that the assertion holds for the final number of primary parastichies
if $\theta$ is rational. A continuity argument implies the 
result in general. 

\begin{figure}[h]\label{figure1}
\epsfysize=5cm
\center{\epsfbox{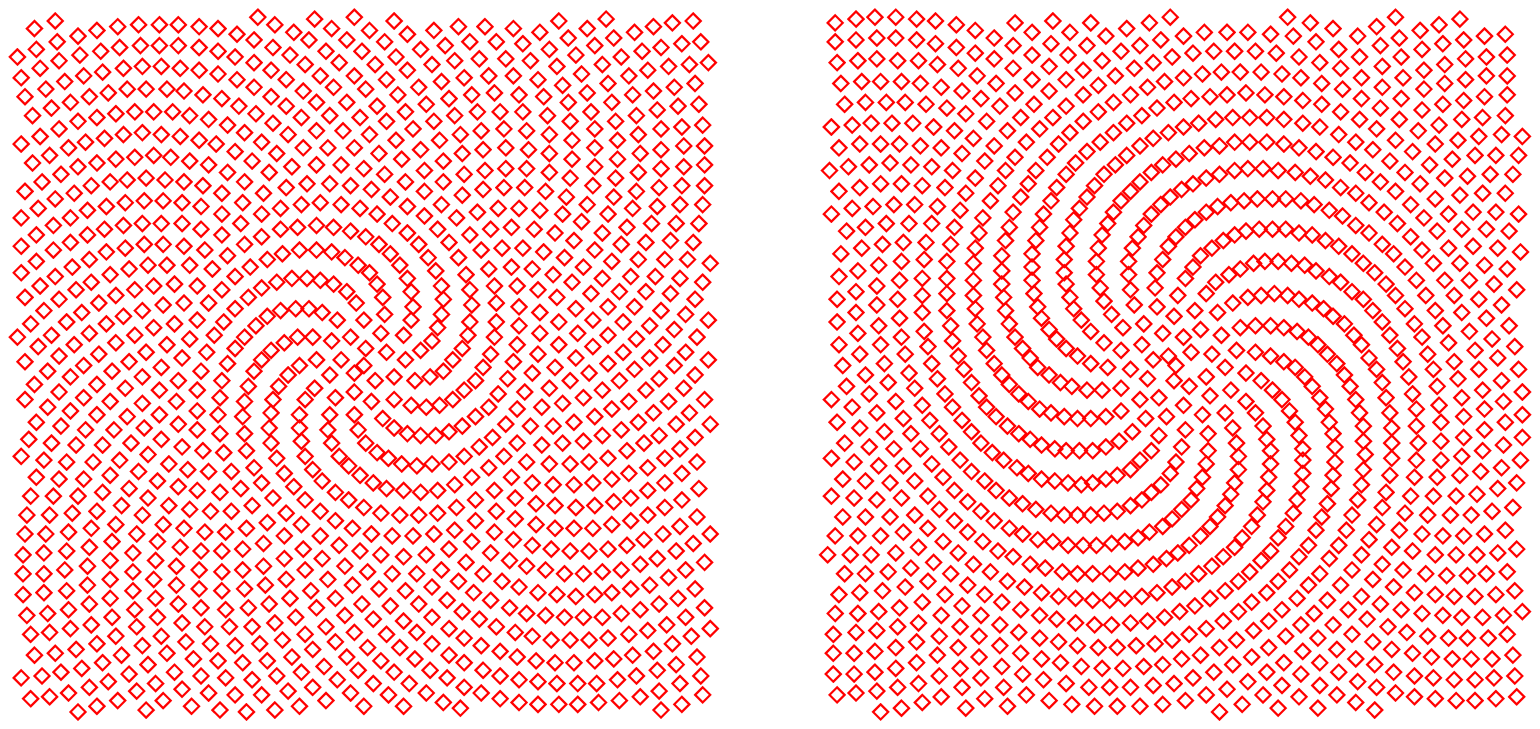}}
\caption{Small points of $\varphi_\theta(\mathbb N)$ for $\theta
=\frac{1765-\sqrt 5}{2858}$ (left side) and $\theta=e^{-1}$ (right side).}
\end{figure}

The left half of Figure 4 shows all points of $\varphi_\theta(\mathbb N)$ 
with real and imaginary parts smaller than 30 for 
$\theta=\frac{1765-\sqrt 5}{2858}$ and $\theta=e^{-1}$.
We have 
$\frac{1765-\sqrt 5}{2858}=[0;1,1,1,1,1,1,3,1,1,1,1,\dots]$
and the first convergents are
$$\frac{0}{1},\frac{1}{1},\frac{1}{2},\frac{2}{3},
\frac{3}{5},\frac{5}{8},\frac{8}{13},\frac{29}{47},\frac{37}{60},\dots\ .$$
The $13$ parastichy spirals corresponding to the denominator $13$
are clearly visible. This is of course due to the \lq\lq large'' coefficient
$3$ in the continued fraction expansion of $\theta$ which leads to three 
parastichy transitions of type II.

Similarly, we get for the parameter 
$\theta=e^{-1}=[0;2,1,2,1,1,4,1,1,6,1,1,8]$ corresponding
to the right half of Figure $4$ the convergents
$$\frac{0}{1},\frac{1}{2},\frac{1}{3},\frac{3}{8},
\frac{4}{11},\frac{7}{19},\frac{32}{87},\frac{39}{106},\dots$$
with a clearly visible parastichy family corresponding to the denominator
$19$. This family is due to the large coefficient $4$ (the family corresponding
to $6$ becomes visible at a larger scale) which causes four
parastichy transitions of type II.

%Example with Lucas sequence $1,3,4,7,11,18,18,29,\dots$.

\subsection{Monodromy}

One can consider two notions of monodromy for phyllotactic sets:

A first notion consists in moving a chosen basis for the local \lq\lq lattice''
by comparing bases of close points in the obvious way. No monodromy
arises in this way: A closed loop gives rise to the identity.

A second, slightly more interesting feature is translational monodromy:
going counterclockwise 
around the origin on a piecewise linear path and stitching 
the obtained lattice elements (with respect to ``bases'' which are 
``moved'' continuously) together, we get at a point 
$a+ib\in \varphi_{\theta}(\mathbb N)$ a vector close to  
$2\pi(-b+ia)$. More precisely, using a basis $V_1,V_2$
associated to primary and secondary 
parastichy spirals, this vector is of the form
$(q_{j-1}+kq_j)V_1\pm q_jV_2$ (with signs depending on the 
sign conventions for $V_1$ and $V_2$) if $\varphi_\theta(\mathbb N)$
contains $q_j$ primary and $(q_{j-1}+kq_j)$ secondary parastichy families
at distance $\sqrt{a^2+b^2}$ from the origin.

%%%%%%%%%%%%%%%%%%%%%%%%%%%%%%%%%%%%%%%%%%%%%%%%%%%%%%%%%%%%%%%%%%%%%%%%%%%%%

\section{Combinatorics of Voronoi-diagrams}\label{sectcombvor}

The material for this Section was motivated by empirical observations
described in \cite{SCR}. This 
Section presents geometric explanations of the observed features
and describes a few related combinatorial facts.

\subsection{Voronoi-diagrams of complex lattices}
Voronoi cells of a complex lattice $\Gamma
=\mathbb Z+\mathbb Z\tau$ corresponding to a complex number $\tau$
in the fundamental domain $\mathcal M$ defined by (\ref{deffunddomM})
yield a tiling of the plane $\mathbb C=\mathbb R+i\mathbb R$ with a
fundamental domain given by the Voronoi cell of the origin.
Voronoi domains are rectangles for 
$\tau=is\in\mathcal M\cap i\mathbb R$ and convex hexagons 
formed by three pairs of parallel opposite edges
otherwise, ie. if the real part of $\tau\in \mathcal M$ is non-zero.

There are thus only two possibilities for the combinatorics 
of the tiling defined by all Voronoi-cells of a 
complex lattice. They are illustrated by the square tiling obtained
by covering the Euclidean plane with unit-squares centered at all
points of $\mathbb Z^2$ and by the honeycomb-tiling consisting of
regular hexagons.

We define the \emph{Voronoi-diagram} of a discrete set $\mathcal S$
in an Euclidean plane $\mathbb E^2$ as the set of all points of $\mathbb E^2$
having at least two closest points in $\mathcal S$ 
at the same minimal distance.
A Voronoi diagram is a plane graph with vertices given by points of 
$\mathbb E^2$  beeing closest to at least three points in $\mathcal S$. 
Edges are given 
by points of $\mathbb E^2$ equidistant to exactly two closest 
elements in $\mathcal S$.
Voronoi cells are open convex polygons (not necessarily bounded) 
with a unique closest point in $\mathcal S$. They define connected 
components in the complement of the Voronoi diagram.

The combinatorics of 
a honeycomb diagram (or more generally of any reasonably locally finite
Voronoi diagram involving only vertices 
of degree $3$) are stable under small perturbations: 
Small independent perturbations of all points of a generic lattice 
corresponding to the honeycomb-case do not change 
the combinatorics of the Voronoi diagram. 
This does not hold for square (or rectangular) lattices: 
a small generic lattice-perturbation
of a lattice $\mathbb Z+is\mathbb Z$ (for $s\geq 1$ a fixed real number)
amounts to splitting all vertices of the grid-graph into two close 
adjacent vertices. For a continuous deformation
$t\longmapsto \mathbb Z+\mathbb Z\tau(t)$
such that $\tau(0)=is$, 
$\Re(\tau(t))<0$, respectively $\Re(\tau(t))>0$, for $t<0$, respectively
for $t>0$, one of the three pairs of opposite 
parallel edges in a Voronoi cell
degenerates to an edge of length $0$, the two other pairs of parallel edges
become orthogonal. Combinatorially, we observe $I-H$ transformations which are
applied simultaneously to all edges in one of the three families of parallel
edges of the Voronoi-tiling, see Figure 5 for an illustration of a continuous
lattice-deformation.
%\ref{figureIH}.
\begin{figure}[h]\label{figureIH}
\epsfysize=4cm
\center{\epsfbox{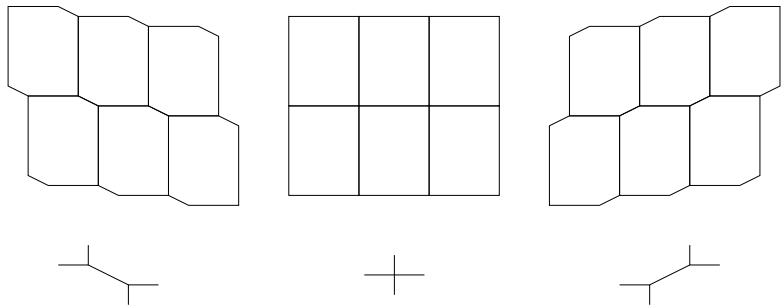}}
\caption{An $I-H$ transformation on Voronoi diagrams of lattices}
\end{figure}
There are thus two possibilities (illustrated by Figure 5)
for \lq\lq resolving'' 
vertices of degree $4$ into pairs of vertices of degree $3$ in a
planar graph. Choosing one of the two possible resolutions
for all vertices of the grid graph (with
vertices $\mathbb Z^2$ and horizontal or vertical edges of unit length)
yields a $3-$regular planar 
graph with cells having at least $4$ and at most $8$ vertices. More precisely,
after separation of points with different resolutions by suitable 
\lq\lq defect-curves'',
the resulting $3-$regular graphs can be achieved as Voronoi-graphs 
of generic perturbations obtained by
small vertical contractions and horizontal 
expansions, respectively horizontal contractions and vertical
expansions, of the connected domains enclosed by defect curves
(for unbounded domains one has to choose
contractions-expansions which are asymptotically very close to the identity).
This construction works of course for any grid-graph 
associated to a partial tiling of a subset of 
$\mathbb R^2$ by identical rectangles. Figure 6 shows an example: defect curves 
are dotted, points of the original grid-graph are replaced by
short fat diagonal
edges with midpoints given by the original lattice points and 
with endpoints joined in the obvious way by (almost) horizontal and vertical 
edges. Observe that the short fat
diagonal edges are of length $O(\epsilon)$ in the Voronoi diagram of a 
generic perturbation of the set $\mathbb Z^2$ moving points less than
$\epsilon$. The remaining set of ``non-short'' or regular edges form a set of 
non-intersecting curves. These piece-wise linear 
curves formed by edges which are alternatingly almost horizontal,
respectively almost vertical, run in ``parallel'' in domains
corresponding to the same type of resolution and are in some sense
``orthogonal'' 
(for an initial grid formed by squares)
in two domains corresponding to different types of resolution. 
They change direction when crossing 
defect curves. Every Voronoi cell contains exactly $4$ segments of 
such curves. Quadrilateral Voronoi cells (arising at suitable
intersections of two defect curves) are enclosed by a unique closed curve 
consisting of four regular edges. Pentagons are almost enclosed by 
such a curve making a half-turn. Octogons are delimited locally by four 
different curves bending away from it.
\begin{figure}[h]\label{figureIH}
\epsfysize=8cm
\center{\epsfbox{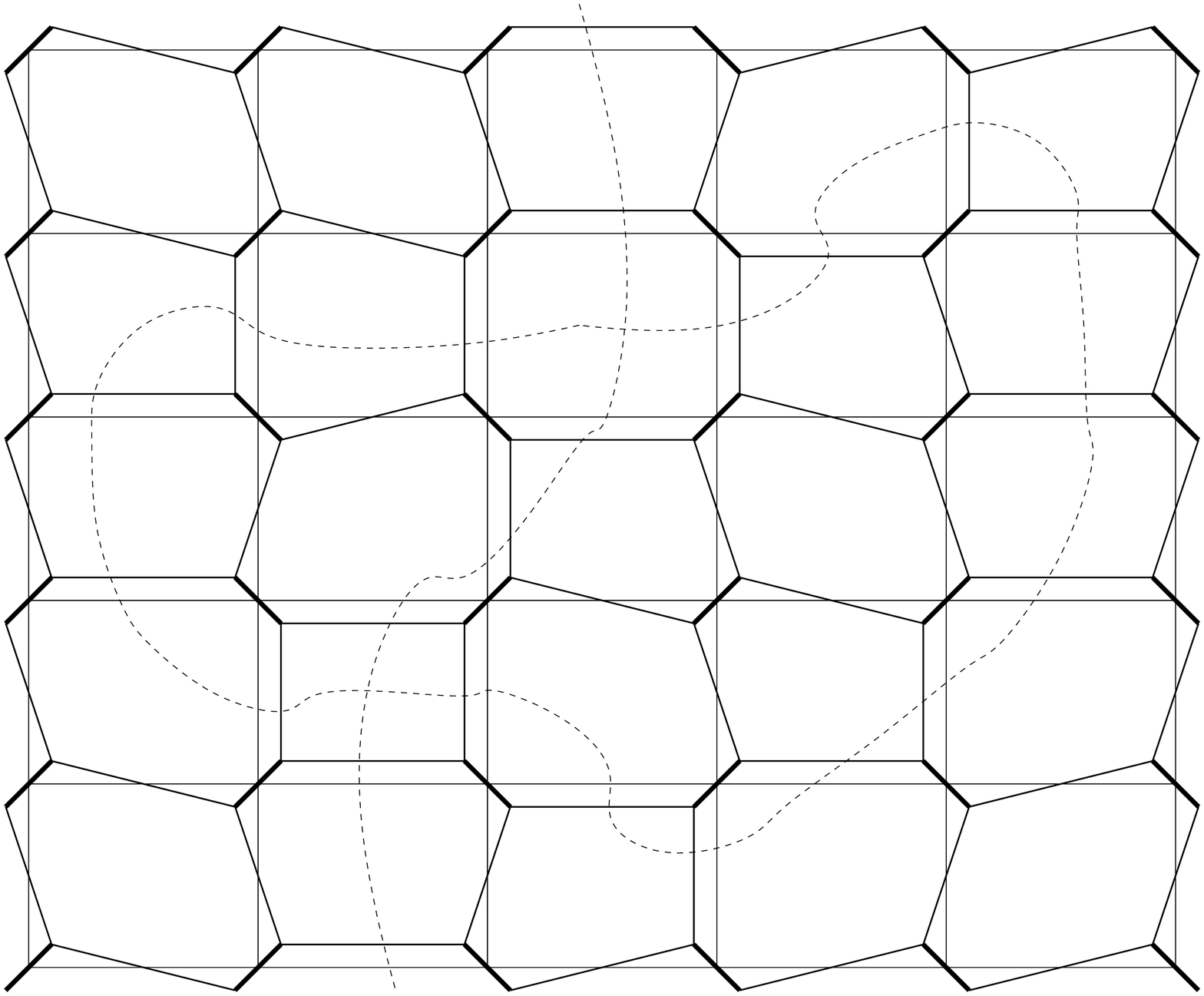}}
\caption{A generic perturbation of a rectangular tiling}
\end{figure}

In the case of a phyllotactic set $\varphi_\theta(\mathbb N)$, 
the defect curves are asympotically given by circles with radii $\sqrt{\rho_i}$
associated to points $\gamma_\theta(\rho_i)$ of the 
phyllotactic geodesic 
corresponding to rectangular lattices (ie. such that $\mathrm{PSL}_2(\mathbb Z)
\left(\gamma_\theta(\rho_i)\right)$ intersects $i\mathbb R$).
Since these circles have asymptotically large radii, defect curves are
locally almost straight lines. Vonoroi cells intersecting a given defect curve 
form a chain with heptagons locally 
nearest to the origin paired
to pentagons at locally maximal distance to the origin. Such ``defect-dipoles''
are joined by chains of ``defect-hexagons'' having two very small sides 
at distance $2$ ``pointing to the next pentagon'' (non defectuous hexagons 
have approximatively parallel opposite sides of roughly equal lengths),
see Figure 7 for a combinatorial illustration not respecting lengths
(small edges would be invisible otherwise) with a straight dashed line 
representing an ideal defect circle of infinite radius. The exterior of
the defect circle is above the dashed line.
\begin{figure}[h]\label{figureIH}
\epsfysize=4cm
\center{\epsfbox{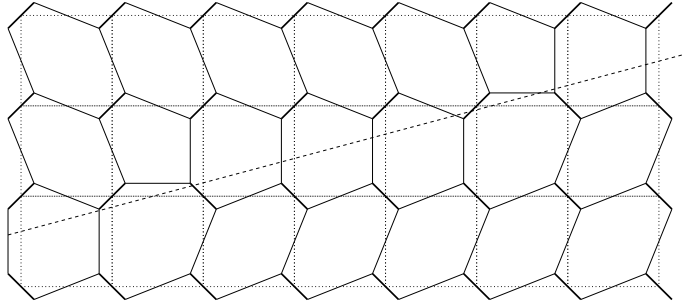}}
\caption{Local combinatorial situation near a defect circle}
\end{figure}

The numbers of defect hexagons and of 
defect dipoles are related in a straightforward way to the numbers
of primary and secondary parastichy spirals. More precisely,
if $\varphi_\theta(\mathbb N)$ has $a$, respectively $b$, parastichy curves
of the two possible types (primary, secondary) with $a<b$ at the defect circle,
then the defect circle meets $a$ defect pentagons and heptagons and $b-a$
defect hexagons. 
Primary and secondary parastichy spirals make up $\frac{2}{3}$ of all edges
(corresponding to all ``vertical'' and ``horizontal'' edges of the unperturbed 
rectangular tiling represented by dotted lines in Figure 7) 
of the dual graph of the Voronoi-graph.
The remaining edges of the Voronoi graph form ``tertiary parastichy spirals''
given by pairs of roughly opposite third-nearest points and running 
\lq\lq in parallel'' to the piecewise-linear 
curves formed by regular edges (and consisting alternatingly of 
\lq\lq horizontal'' and \lq\lq vertical edges'' in Figure 7). 

Figure 1 of \cite{SCR} is a picture in the case where $\theta$ is 
the golden mean. Two consecutive 
parastichy transitions of type II are always separated by a defect circle.
Two parastichy transitions of type II separated by a
parastichy transition of type I, are separated by either
one or no defect circle.
In the case $\theta=\frac{1+\sqrt{5}}{2}$ there is (asymptotically)
always such a defect circle 
(which coincides almost with a parastichy transformation of type I).
Defect-dipoles are separated at most by a unique defect-hexagon.

%\begin{figure}[h]\label{figureIH}
%\epsfysize=14cm
%\epsfysize=15cm
%\center{\epsfbox{voron130.eps}}
%\center{\epsfbox{resolution1.eps}}
%\caption{Local combinatorial situation near a defect circle}
%\end{figure}

%The point $n=2280$ corresponds to a point of $\mathcal M$ very close to $i$.
%$n=2277$ corresponds to a good point of the phyllotactic set.
%(five outside, seven inside).
%%%%%%%%%%%%%%%%%%%%%%%%%%%%%%%%%%%%%%%%%%%%%%%%%%%%%%%%%%%%%%%%%%%%%%%%%%%%%

\section{Chromatic properties}\label{sectchrom}

A \emph{Tait graph} is a $3-$regular graph with $3-$coloured edges 
such that three edges of all three edge-colours meet at every vertex.
Plane Tait graphs are essentially the same as $4-$coloured generic maps 
using the following classical trick. Identify the $3$ edge-colours with 
the non-zero elements of Klein's Viergruppe $\mathbb V=
\mathbb Z/2\mathbb Z\times
\mathbb Z/2\mathbb Z$. Up to a colour-permutation, there is then a unique
colouring with colours $\mathbb V$ such that adjacent regions coloured 
$\alpha$ and $\beta$ are separated by an edge of colour $\alpha+\beta$.

Colouring parallel edges of the Voronoi-diagram of a generic lattice 
with three colours yields a Tait graph. The associated $4-$colouring
of the Voronoi-diagram has the property that the boundary of
every hexagon meets 
all three other colours cyclically. 
Opposite parallel edges separate thus a given hexagon from 
two neighbours sharing the same colour.
This colouring is unique up to colour-permutation. 
(Observe that the Voronoi diagram of the hexagonal lattice has also
a nice essentially unique colouring involving only $3$ colours
associated to a suitable morphism from the hexagonal lattice into 
a cyclic group of order $3$).

This Tait colouring exists of course for a generic radius of the phyllotactic 
set and can be extended uniquely to a Tait colouring of the interior
except at the boundary of the Voronoi-domain containing the origin.
(Indeed, the situation at defect circles is quite easy to understand:
there are ``vertical'', ``horizontal'' and very short (asymptotically 
infinitesimal) ``defect'' edges defining a coherent Tait-colouring
which is unique up to colour-permutations.)

The Voronoi-domain of the phyllotactic set has thus a locally canonical
$4-$colouring (except at the origin) displaying chromatic monodromy 
when trying to extend this canonical colouring along a loop 
around the origin: Returning to the starting point we get a final colouring 
which is different. More precisely, the four colours get exchanged by pairs.

\begin{figure}[h]\label{figureIH}
\epsfysize=8cm
%\epsfysize=15cm
%\center{\epsfbox{voron130.eps}}
\center{\epsfbox{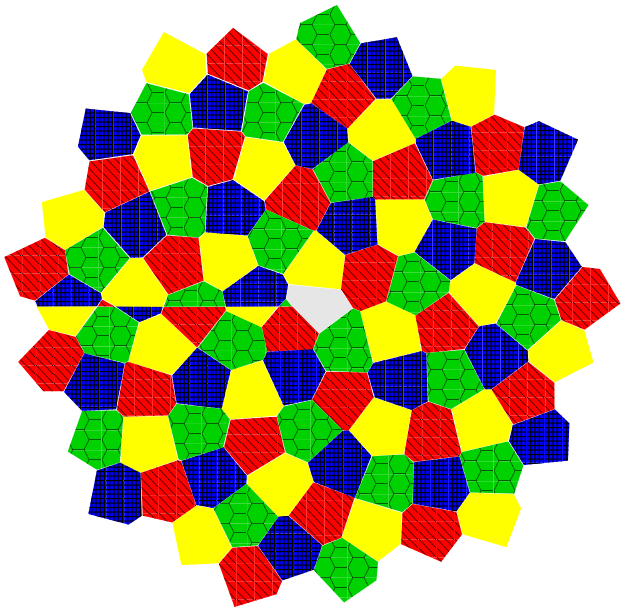}}
\caption{$4-$colouring of the first Voronoi-cells for $\theta=\frac{1+\sqrt{5}}{2}$.}
\end{figure}

Figure 8 shows the canonical colouring with chromatic-monodromy
concentrated at the real negative halfline.

A slightly different (but in fact equivalent) construction of this
canonical local $4-$colouring is as follows: Voronoi-cells of a lattice 
$\Lambda\subset \mathbb C$ are always $4-$coloured using the homomorphism
$\Lambda\longrightarrow \Lambda/2\Lambda\sim \mathbb V$. Affine lattice
transformations preserving $\Lambda$ induce affine transformations
of the colour group $\mathbb V$. More precisely, a translation of 
$\Lambda$ by an element $\lambda\not\in 2\Lambda$ induces a translation of all 
colours in $\mathbb V\sim\Lambda/2\Lambda$
by the corresponding element $\lambda\pmod{2\Lambda}$.
We get thus asymptotically a local 
$4-$colouring of the phyllotactic set which is canonical up to addition of
a constant in the colour group $\mathbb V$.
Pushing the affine colour translations along a loop around the origin
leads to the monodromy-translation of $\mathbb V$ partitioning 
$\mathcal V$ into two orbits of pairwise exchanged colours.

Colouring black and white the two monodromy-orbits of $\mathbb V$,
we get a partition into black and white cells 
of all Voronoi-domains not containing the origin. Two adjacent black
cells or two adjacent white cells are always separated by edges of the 
same ``monodromy'' colour associated to the monodromy translation of the
colour-group $\mathbb V$.
\begin{figure}[h]\label{figureIH}
\epsfysize=8cm
%\epsfysize=15cm
%\center{\epsfbox{voron130.eps}}
\center{\epsfbox{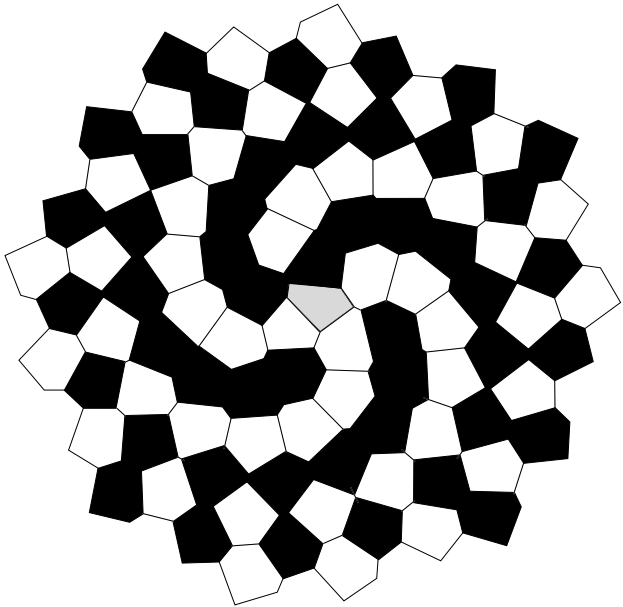}}
\caption{The black-white colouring obtained after identification of 
pairs of colours exchanged by monodromy.}
\end{figure}
Figure 9 displays the associated black and white cells. Centers (ie.
elements of the phyllotactic set defining the Voronoi-cells) of white
cells correspond to images under the phyllotactic map $\varphi_\theta$
(with $\theta=\frac{1+\sqrt{5}}{2}$) of odd natural numbers. Centers of black
cells correspond to images of non-zero even natural numbers. 
Black and white cells form stripes which degenerate
at every third defect circles into a checkerboard-like situation.
Indeed, defect circles are alternatingly associated to one of the three colours
(colouring the asymptotically infinitesimal edges near the defect circle).
At defect circles corresponding to the monodromy colour (ie. where the edges
of the monodromy colour degenerate into very short edges), the picture
degenerates into a checkerboard. Elsewhere, we get alternating black and white
``stripes'' running (more or less) in parallel to a parastichy family.
The total number of stripes (of both colours) is always an even 
Fibonacci number. The lower half of Figure 8 displays
the situation for the even Fibonacci number 8 (the situation for 
the next even Fibonacci 
number, $34$ is just visible at the boundary of the image).
Since indices of consecutive even Fibonacci numbers in the Fibonacci sequence
differ by $3$ which is odd, two
consecutive ``stripy'' regions rotate visually in opposite directions.

Choosing other pairwise identifications (orbits of $\mathbb V$
under addition of one of the other two non-zero elements of $\mathbb V$)
of the four colours represented by $\mathbb V$ leads locally to stripy regions
for other parastichy families. Every parastichy family can be represented 
locally in this way but there is of course a global error, accounted
for by chromatic monodromy, if the number
of curves in the family is odd. The birth-point and the dying point, 
corresponding to radially consecutive checkerboardy
regions of such a family, are at every third defect circle.
They correspond to the moment where the centers of adjacent black (or white)
cells no longer correspond to pairs of third-nearest lattice points
in the linearized lattice. They are thus birth-places (or grave-yards) for
``tertiary'' parastichy families which transform later during a type
II transition into secondary parastichies.

The chromatic constructions described above work of course for many other real 
divergence angles: it is enough that the corresponding real number 
has a continued fraction expansion with coefficients not growing too
fast. This ensures the existence of a (essentially unique)
Tait colouring far from the origin. Near the origin, this Tait
colouring breaks down. It gives however rise asymptotically 
to a local $4-$coulouring (coming from the morphism 
$\Lambda\longrightarrow \Lambda/2\Lambda\sim \mathbb V$)
of Voronoi-cells which displays in general a non-trivial colour-monodromy
around the origin given by a monodromy translation (exchanging
colours pairwise) in $\mathbb V$.

\begin{rem} Whorled phyllotaxis (see Section \ref{subswhorled}) 
of even order leads (for suitable divergence angles) 
to asymptotic Tait colourings with trivial colour-monodromy.
\end{rem}

\begin{rem} Choosing a cyclic order on the edge-colours of a finite 
Tait graph yields a compact surface by gluing discs on all closed 
oriented paths defined by edges with the prescribed cyclical colour-scheme. 

Defect lines give rise to genus proportional to a suitably defined 
combinatorial length.
This gives in particular a topological interpretation of the $2-$dimensional
Ising model with periodic boundary conditions: The associated Hamiltonian
(in the case of no external magnetic field) is
proportional to the genus of the associated
compact surface obtained by ``opening'' all vertices according to their sign.
This leads to a canonical Tait graph allowing the construction of 
a compact surface.
\end{rem}

%%%%%%%%%%%%%%%%%%%%%%%%%%%%%%%%%%%%%%%%%%%%%%%%%%%%%%%%%%%%%%%%%%%%%%%%%%%%%

\section{Other models}\label{sectothermod}

\subsection{Cylindric and logarithmic models}

In \cite{Cox72} Coxeter, following \cite{BB}, 
models the structure of pineapples or pine-cones 
by approximating their shape with a cylinder which he develops on
the plane thus getting an infinite strip of a lattice. The visible
features (scales) on pineapples are
the Voronoi cells of this lattice. The obtained lattice should
stay close to the hexagonal lattice which has optimal packing and covering
properties.
Working with a cylinder of circumference $2\pi$, we have thus to choose
the optimal divergence angle $2\pi\theta$ such that the complex lattice
$\mathbb Z+\mathbb Z(\theta+\epsilon i)$ is close to the hexagonal
lattice for small $\epsilon$. There is no exact control over the value 
of $\epsilon$ since pineapples or pine cones are not exact cylinders.
We should thus choose the value of $\theta$ such that the geodesic
$t\longmapsto \theta+ti$ is overall optimal for small positive $\epsilon$.
For a fixed value of $\theta$, the map $\epsilon\longmapsto 
\theta+\epsilon i$ defines again a geodesic of the hyperbolic 
half-plane.

The best choice is of course again
given by $\lambda=\frac{1+\sqrt{5}}{2}2\pi$ (or by its negative) modulo $2\pi$
yielding a geodesic which is asymptotically close to the shortest
periodic geodesic 
$$(0,\infty)\ni t\longmapsto \frac{1}{2}
\frac{(1+it)+(-1+it)\sqrt 5}{1+it}$$
of the modular curve $\mathrm{PSL}_2(\mathbb Z)\backslash\mathbb H$.

\begin{rem}
Identifying an infinitely long 
cylinder of circumference $1$ with the quotient space of 
$\mathbb C$ under translations of the form $2i\pi \mathbb Z$, 
the usual exponential function transforms the cylindric model into 
the logarithmic model with points on a logarithmic ontogenetic spiral
defined by $\mathbb N\ni n\longmapsto\rho^ne^{2i\pi \theta n}$.
\end{rem}

\subsection{van Iterson's disc-packing model}

In \cite{It} 
van Iterson considers periodic packings of equal discs on cylinders
such that every disc touches (at least) two pairs of adjacent discs.

Levitov in \cite{Lev} observes that 
the associated lattices correspond to elements of norm $1$ in
the modular domain and that the corresponding Teichm\"uller space 
(given by equivalence classes of lattices endowed with a
positively-oriented basis) is a $3-$regular tree with 
mid-edges given by the 
$\mathrm{PSL}_2(\mathbb Z)-$orbit of $i$ (corresponding to the square-lattice)
and vertices given by 
the $\mathrm{PSL}_2(\mathbb Z)-$orbit of $\frac{1+i\sqrt{3}}{2}$
(corresponding to the hexagonal lattice).
The rooted subtree defined in the quarter-plane of $\mathbb C$
defined by all elements of $\mathbb H$ with non-negative real parts
is then in natural bijection with the so-called Farey-tree.
The optimal approximatively straight choice for a path on this rooted tree
corresponds of course again to the golden mean and is given by alternating
left- and right-turns at every bifurcation, as also observed by Levitov 
who gives a physical explanation based on energy levels of this fact.

\begin{rem} van Iterson's model involves sphere packings which are 
\lq\lq locally optimal'' in the sense that every disc of the packing (almost)
touches four other discs. The associated Voronoi cells have however not 
asymptotically equal are\ae.

The model determined by the phyllotactic map $\varphi_\theta$
gives Voronoi cells with asymyptotically equal are\ae
but leads to discs in the corresponding sphere packing
which are almost all isolated. 
Exceptions are occuring at parastichy transitions of type $I$
(happening asymptotically at the square lattice if $\theta$ is
the golden ratio).
\end{rem}

\subsection{Opposite and whorled phyllotaxis}\label{subswhorled}

The content of this paper can be applied in a straightforward way 
to a whorled model of the sunflower: We consider a natural number
$d\geq 2$ (the case $d=2$ is called \lq\lq opposite phyllotaxis'',
the case $d>3$ \lq\lq whorled phyllotaxis'' in botanics) and we denote by 
$\mathbb U_d$ the set of all $d$ complex $d-$th roots of $1$
given by the zeros of the polynomial $z^d-1$. 
The set $\mathbb U_d\varphi_\theta\left(\frac{1}{d}\mathbb N\right)$
defined by the map
$$\mathbb U_d\times \mathbb N\ni(\omega,n)\longmapsto
\omega \sqrt{\frac{n}{d}}e^{2i\pi\theta n/d}$$
(invariant under the obvious isometrical action defined by 
the multiplicative subgroup $\mathbb U_d$ of $\mathbb C^*$) 
is then locally essentially a rescaling by a factor $\frac{1}{d}$ 
(followed by a rotation) of the 
phyllotactic set $\varphi_\theta(\mathbb N)$. 
Up to reparametrization
(replacing $t$ with $td$),
it admits the same
phyllotactic geodesic. Numbers of parastichy curves in a given family 
are simply multiplied by $d$ with respect to the corresponding number 
in $\varphi_\theta(\mathbb N)$.
%%%%%%%%%%%%%%%%%%%%%%%%%%%%%%%%%%%%%%%%%%%%%%%%%%%%%%%%%%%%%%%%%%%%%%%%%

\section{Testing the existence of phyllotactic geodesics in real 
sunflowers}\label{secttest}

Phyllotactic geodesic are perhaps a mere mathematical artefact
due to the use of the model maps $\varphi_\theta$. 
This Section sketches a test probing the 
reality of the theory.

A first step is of course gathering real data, consisting of a fair number 
of pictures of large flawless sunflower-capitula. These pictures should
be enriched by adding as smoothly as possible 
(using perhaps splines or trigonometric functions and a least square 
method) all visible parastichy spirals.
Intersections of transversal parastichy spirals should now 
be taken as the centers of seeds. Points near the center
can be neglected.

We can check adequacy of $\varphi_\theta(\mathbb N)$ for sunflowers
as follows: Determine for each picture
(endowed with a complex coordinate system)
parameters $A\in \mathbb C,C\in\mathbb C^*,\theta,\gamma\in \mathbb R$ giving 
the best least square 
approximation of the obtained seed-centers with a suitable set of 
points of the form
$$\mathbb N\ni n\longmapsto A+C\sqrt{n+\gamma}\ e^{2i\pi\theta n},$$
supposing that the pictures have no distorsions 
(additional parameters are necessary otherwise).
If this approximation is nearly perfect, the sunflower map
is an accurate description of reality and the existence of 
phyllotactic geodesics is confirmed.
A failure or a bad match does however not contradict the existence 
of phyllotactic geodesics but forces us 
to compute points of \lq\lq hypothetical'' geodesics using 
the real data-sets instead of the 
model set $\varphi_\theta(\mathbb N)$.

This can be achieved as follows: For each point $P$ neither on the boundary
nor in the center of the sunflower,
we determine pairs of points $a,A$ and $b,B$ adjacent to
$P$ with $a,A$ on one parastichy spiral through $P$
and $b,B$ on the other, transversal parastichy spiral through $P$.
The linearized lattice at $P$ is then approximatively
given by $\mathbb Z\frac{A-a}{2}+\mathbb Z
\frac{B-b}{2}$. This allows the computation of the corresponding 
modular invariant by considering the point
of the modular curve represented by $\pm \frac{A-a}{B-b}$ (for the
unique sign choice leading to a strictly positive imaginary part). 
Suitable lifts of these points to $\mathbb H$
should now lie close to a hyperbolic geodesic which can be guessed by 
least square approximation.

\hskip.5cm

{\bf Acknowledgements} I would like to thank David Speyer who 
started my interest in phyllotaxis by proposing 
$\varphi_\theta(\mathbb N)$ as an interesting configuration 
to consider in relation with Question 3307 of Mathoverflow, see \cite{M}, 
to Tanguy 
Rivoal who prompted me to write up the details and to Pierre de la Harpe
for useful remarks.

%%%%%%%%%%%%%%%%%%%%%%%%%%%%%%%%%%%%%%%%%%%%%%%%%%%%%%%%%%%%%%%%%%%%%%%%%%%%%

%http://www.math.smith.edu/phyllo/    contains interesting material

\noindent Roland BACHER, Universit\'e Grenoble I, CNRS UMR 5582, Institut 
Fourier, 100 rue des maths, BP 74, F-38402 St. Martin d'H\`eres, France.

\noindent e-mail: Roland.Bacher@ujf-grenoble.fr

\end{document}